\documentclass[a4paper,11pt]{amsart}

\usepackage{amsmath}
\usepackage{amssymb}
\usepackage{amsthm}
\usepackage[abbrev,alphabetic]{amsrefs}
\usepackage{amscd}
\usepackage[dvipdfmx]{graphicx}
\usepackage{hyperref}

\usepackage{amsthm}
\usepackage{color}
\usepackage[normalem]{ulem} 
\newcommand\redout{\bgroup\markoverwith
{\textcolor{red}{\rule[.5ex]{2pt}{0.4pt}}}\ULon}
\newcommand\blueout{\bgroup\markoverwith
{\textcolor{blue}{\rule[.5ex]{2pt}{0.4pt}}}\ULon}

%; Secondary }


\newtheorem{thm}{Theorem}[section]
\newtheorem{lem}[thm]{Lemma}

\newtheorem{cor}[thm]{Corollary}
\newtheorem{prop}[thm]{Proposition}

\newtheorem{step}{Step}

\theoremstyle{definition}

\newtheorem{dfn}[thm]{Definition}
\newtheorem{nota}[thm]{Notation}

\theoremstyle{remark}

\newtheorem{rem}[thm]{Remark}
\newtheorem*{ack}{Acknowledgements}

\newcommand{\MO}{\mathcal{O}}

\newcommand{\Q}{\mathbb{Q}}
\newcommand{\R}{\mathbb{R}}
\newcommand{\Z}{\mathbb{Z}}

\newcommand{\red}[0]{{\operatorname{red}}}
\newcommand{\Ex}[0]{{\operatorname{Ex}}}

\DeclareMathOperator{\Supp}{Supp}

\DeclareMathOperator{\Spec}{Spec}

\DeclareMathOperator{\Pic}{Pic}

\DeclareMathOperator{\NE}{NE}

\DeclareMathOperator{\Mov}{Mov}
\DeclareMathOperator{\coeff}{coeff}

\DeclareMathOperator{\Nklt}{Nklt}

\makeatletter
  
  \@addtoreset{equation}{thm}
  \makeatother

\title[Minimal model program for log canonical threefolds]
{Minimal model program for log canonical threefolds in positive characteristic}

\author{Kenta Hashizume}
\address{Department of Mathematics, Graduate School of Science, 
Kyoto University, Kyoto 606-8502, Japan}
\email{hkenta@math.kyoto-u.ac.jp}

\author{Yusuke Nakamura}
\address{Graduate School of Mathematical Sciences, 
the University of Tokyo, 3-8-1 Komaba, Meguro-ku, Tokyo 153-8914, Japan.}
\email{nakamura@ms.u-tokyo.ac.jp}

\author{Hiromu Tanaka}
\address{Graduate School of Mathematical Sciences, 
the University of Tokyo, 3-8-1 Komaba, Meguro-ku, Tokyo 153-8914, Japan.} 
\email{tanaka@ms.u-tokyo.ac.jp}

\begin{document}

\maketitle

\begin{abstract}
Given a three-dimensional projective log canonical pair over a perfect field of characteristic larger than five, 
there exists a minimal model program that terminates after finitely many steps. 
\end{abstract}

\tableofcontents

\section{Introduction}

Recently there have been large developments on the three-dimensional minimal model program of characteristic $p>5$. 
The first remarkable achievement was made by Hacon and Xu. 
They proved the existence of minimal models for terminal threefolds over 
an algebraically closed field $k$ of characteristic $p>5$ \cite{HX15}. 
Actually, they succeeded in establishing the so-called generalised minimal model program 
for terminal threefolds $X$ over $k$ such that $K_X$ is pseudo-effective. 
Then Cascini, Xu and the third author dropped the assumption that 
$K_X$ is pseudo-effective \cite{CTX15}. As a consequence, 
it turned out that an arbitrary terminal threefold over $k$ is birational to either a minimal model or a Mori fibre space. 
Finally, Birkar and Waldron established the minimal model program for klt threefolds over $k$ \cite{Bir16}, \cite{BW17}.

The main objective of this paper is to give a generalisation of 
the three-dimensional minimal model program over $k$ to the log canonical case. 
Furthermore, we treat perfect base fields, the relative setting and 
the non-$\Q$-factorial case. 
More specifically, our main theorem is as follows.

\begin{thm}[cf.\ Theorem \ref{t-lc-mmp}]\label{intro-mmp}
Let $k$ be a perfect field of characteristic $p>5$. 
Let $(X, \Delta)$ be a three-dimensional log canonical pair over $k$, 
where $\Delta$ is an $\R$-divisor. 
Let $f:X \to Z$ be a projective $k$-morphism to a quasi-projective $k$-scheme $Z$. 
Then there exists a $(K_X+\Delta)$-MMP over $Z$ that terminates. 

In particular, if $X$ is $\Q$-factorial, then 
there is a sequence of birational maps of three-dimensional normal varieties:  
\[
X=:X_0 \overset{\varphi_0}{\dashrightarrow} X_1 \overset{\varphi_1}{\dashrightarrow} \cdots \overset{\varphi_{\ell-1}}{\dashrightarrow} X_{\ell}
\]
such that if $\Delta_i$ denotes the proper transform of $\Delta$ on $X_i$, then 
the following properties hold:  
\begin{enumerate}
\item 
For any $i \in \{0, \ldots, \ell\}$, 
$(X_i, \Delta_i)$ is a $\Q$-factorial log canonical pair which is projective over $Z$.
\item 
For any $i \in \{0, \ldots, \ell-1\}$, 
$\varphi_i:X_i \dashrightarrow X_{i+1}$ is either a $(K_{X_i}+\Delta_i)$-divisorial contraction over $Z$ or a $(K_{X_i}+\Delta_i)$-flip over $Z$. 
\item 
If $K_X+\Delta$ is pseudo-effective over $Z$, then $K_{X_{\ell}}+\Delta_{\ell}$ is nef over $Z$. 
\item 
If $K_X+\Delta$ is not pseudo-effective over $Z$, then 
there exists a $(K_{X_{\ell}}+\Delta_{\ell})$-Mori fibre space $X_{\ell} \to Y$ over $Z$. 
\end{enumerate}
\end{thm}

\begin{rem}\label{r-intro-mmp}
Theorem \ref{intro-mmp} is known for the following cases: 
\begin{enumerate}
\item 
$(X, \Delta)$ is $\Q$-factorial klt and $\Delta$ is a $\Q$-divisor \cite[Theorem 2.13]{GNT}. 
\item 
$k$ is algebraically closed, $Z$ is projective over $k$, and $K_X+\Delta \equiv_f D$ 
for some effective $\R$-Cartier $\R$-divisor $D$ \cite[Corollary 1.8]{Wal}. 
\end{enumerate}
\end{rem}

To prove the main theorem (Theorem \ref{intro-mmp}), 
we also establish some fundamental results such as the cone theorem and the base point free theorem under the same generality. 

\begin{thm}[cf.\ Theorem \ref{t-lc-cone}, Theorem \ref{t-lc-cone2}]\label{intro-cone}
Let $k$ be a perfect field of characteristic $p>5$. 
Let $(X, \Delta)$ be a three-dimensional log canonical pair over $k$, 
where $\Delta$ is an $\R$-divisor. 
Let $f:X \to Z$ be a projective $k$-morphism to a quasi-projective $k$-scheme $Z$. 
Then there exists a countable set $\{C_i\}_{i \in I}$ of curves on $X$ 
which satisfies the following properties: 
\begin{enumerate}
\item 
$f(C_i)$ is a point for any $i \in I$. 
\item 
$\overline{\NE}(X/Z)=\overline{\NE}(X/Z)_{K_X+\Delta \ge 0}+\sum_{i \in I} \R_{\geq 0}[C_i].$
\item 
If $A$ is an $f$-ample $\R$-Cartier $\R$-divisor on $X$, then 
there exists a finite subset $J$ of $I$ such that 
\[
\overline{\NE}(X/Z)=\overline{\NE}(X/Z)_{K_X+\Delta+A \ge 0}+\sum_{j \in J} \R_{\geq 0}[C_j].
\]
\item 
If $k$ is algebraically closed, then the set $\{C_i\}_{i \in I}$ can be chosen so that 
$C_i$ is a rational curve such that $0<-(K_X+\Delta) \cdot C_i \leq 6$ for any $i \in I$. 
\end{enumerate}
\end{thm}

\begin{thm}[cf.\ Theorem \ref{t-lc-big-bpf2}, Theorem \ref{t-R-rel-bpf1}, Theorem \ref{t-R-rel-bpf2}]\label{intro-bpf}
Let $k$ be a perfect field of characteristic $p>5$. 
Let $(X, \Delta)$ be a three-dimensional log canonical pair over $k$, 
where $\Delta$ is an $\R$-divisor. 
Let $f:X \to Z$ be a projective $k$-morphism to a quasi-projective $k$-scheme $Z$. 
Let $L$ be an $f$-nef $\R$-Cartier $\R$-divisor on $X$. 
Assume that one of the following conditions holds: 
\begin{enumerate}
\item $L-(K_X+\Delta)$ is $f$-semi-ample and $f$-big. 
\item $L=K_X+\Delta$ and $L$ is $f$-big. 
\item $L=K_X+\Delta$ and $\Delta$ is an $f$-big $\R$-Cartier $\R$-divisor. 
\end{enumerate}
Then $L$ is $f$-semi-ample. 
\end{thm}

\begin{rem}
Theorem \ref{intro-cone} and Theorem \ref{intro-bpf}(1)(2) are known 
for the case when $k$ is algebraically closed and $Z$ is projective over $k$ 
\cite[Theorem 1.1, Theorem 1.2]{Wal}. 
\end{rem}

Although the main theorem (Theorem \ref{intro-mmp}) only asserts the existence 
of a minimal model program that terminates, 
we also obtain several results on termination of flips.

\begin{thm}[cf.\ Theorem \ref{t-lc-eff-mmp2}, Theorem \ref{t-pseff-lc-mmp}]\label{intro-termination}
Let $k$ be a perfect field of characteristic $p>5$. 
Let $(X, \Delta)$ be a three-dimensional log canonical pair over $k$, 
where $\Delta$ is an $\R$-divisor. 
Let $f:X \to Z$ be a projective $k$-morphism to a quasi-projective $k$-scheme $Z$. 
Then the following hold. 
\begin{enumerate}
\item 
If $K_X+\Delta \equiv_f D$ for some effective $\R$-Cartier $\R$-divisor $D$ on $X$, 
then there exists no infinite sequence that is a $(K_X+\Delta)$-MMP over $Z$. 
\item 
Let $A$ be an ample $\R$-Cartier $\R$-divisor. 
If $K_X+\Delta$ is $f$-pseudo-effective,  
then there exists no infinite sequence that is a $(K_X+\Delta)$-MMP over $Z$ with scaling of $A$. 
\end{enumerate}
\end{thm}

%\subsection{\textcolor{red}{Related results}} 

We now summarise some known results related to this paper. 
Several results toward the abundance conjecture for threefolds has been established 
(\cite{DW}, \cite{Wal17}, \cite{XZ}, \cite{Zha}). 
For threefolds of characteristic $p>5$, 
the Iitaka conjecture has been proven to be true if the generic fibre is smooth \cite{EZ}. 
For three-dimensional log canonical pairs of characteristic $p>5$, Das and Hacon prove that 
minimal lc centres are normal \cite{DH16}.

\begin{ack}
The authors would like to thank Caucher Birkar, Paolo Cascini, Osamu Fujino, 
Yoshinori Gongyo, and Joe Waldron for many useful discussions and comments. 
The authors also thank the referee for reading the manuscript carefully and for suggesting several improvements. 
The first author was supported by JSPS KAKENHI Grant Number JP16J05875. 
The second author was supported by JSPS KAKENHI Grant Number 18K13384. 
The third author was funded by EPSRC and JSPS KAKENHI  Grant Number 18K13386. 
\end{ack}

\section{Preliminaries}

\subsection{Notation}\label{ss-notation}

In this subsection, we summarise notation used in this paper. 

\begin{itemize}
\item We will freely use the notation and terminology in \cite{Har77} 
and \cite{Kol13}. 

\item For a scheme $X$, its {\em reduced structure} $X_{\red}$ 
is the reduced closed subscheme of $X$ such that the induced morphism 
$X_{\red} \to X$ is surjective.

\item For an integral scheme $X$, 
we define the {\em function field} $K(X)$ of $X$ 
to be $\MO_{X, \xi}$ for the generic point $\xi$ of $X$. 

\item For a field $k$, 
we say that $X$ is a {\em variety over} $k$ or a $k$-{\em variety} if 
$X$ is an integral scheme that is separated and of finite type over $k$. 
We say that $X$ is a {\em curve} over $k$ or a $k$-{\em curve} 
(resp.\ a {\em surface} over $k$ or a $k$-{\em surface}, 
resp.\ a {\em threefold} over $k$) 
if $X$ is an $k$-variety of dimension one (resp.\ two, resp.\ three). 

\item Let $f:X \to Z$ be a projective morphism of noetherian separated schemes, 
where $X$ is an integral normal scheme. 
For a Cartier divisor $L$ on $X$, 
we say that $L$ is $f$-{\em free} 
if the induced homomorphism $f^*f_*\MO_X(L) \to \MO_X(L)$ is surjective. 
Let $M$ be an $\R$-Cartier $\R$-divisor on $X$. 
We say that $M$ is $f$-{\em ample} (resp.\ \textit{$f$-semi-ample}) 
if we can write $M=\sum_{i=1}^r a_iM_i$ 
for some $r \geq 1$, positive real numbers $a_i$ and 
$f$-ample (resp.\ $f$-free) Cartier divisors $M_i$. 
We say that $M$ is $f$-{\em big} if we can write $M=A+E$ 
for some $f$-ample $\R$-Cartier $\R$-divisor $A$ and effective $\R$-divisor $E$. 
We define $f$-nef $\R$-divisors in the same way as in \cite[Definition 1.4]{Kol13}. 
We say that $M$ is \textit{$f$-numerically trivial}, denoted by $M \equiv_f 0$ or $M\equiv_{Z} 0$, 
if both $M$ and $-M$ are $f$-nef.

Given $\R$-Cartier $\R$-divisors $D_1$ and $D_2$ on $X$ such that 
$D_1+ \lambda_0D_2$ is $f$-nef for some $\lambda_0 \in \R_{>0}$, 
the {\em $f$-nef threshold} of $(D_1, D_2)$ is the non-negative real number $\lambda$ defined by 
\[
\lambda:=\inf\{\mu \in \R_{\geq 0}\,|\,D_1+\mu D_2\text{ is }f\text{-nef}\}.
\]
It is well-known that $\inf$ can be replaced by $\min$ in this definition 
if $Z$ is of finite type over a field. 
\item 
Let $\Delta$ be an $\R$-divisor on an integral normal separated noetherian scheme 
and let $\Delta=\sum_{i \in I} r_i D_i$ be its irreducible decomposition. 
We define $\Delta ^{\ge 1} := \sum_{i \in I, r_i \ge 1} r_i  D_i$ and $\Delta ^{\wedge 1} := \sum_{i \in I} r_i ' D_i$ where $r_i ' := \min \{ r_i, 1  \}$. 
We also define $\Delta ^{> 1}$ and $\Delta ^{< 1}$ similarly. 
Furthermore, we set $\{ \Delta \} := \Delta - \lfloor \Delta \rfloor$. 
\item 
Let $X$ be a normal variety $X$ over a field and 
let $X'$ be a non-empty open subset $X'$ of $X$. 
For a prime divisor $F$ on $X'$, 
the {\em closure} $\overline F$ of $F$ in $X$ is the prime divisor on $X$ 
whose generic point is equal to the generic point of $F$. 
For an $\R$-divisor $D$ on $X'$, 
the {\em closure} $\overline D$ of $D$ in $X$ is defined as $\sum_{i \in I} r_i \overline D_i$, where $D=\sum_{i \in I} r_i D_i$ is the irreducible decomposition of $D$. 
\item 
Let $k$ be a field and let $f:X \to Y$ be a birational $k$-morphism 
of normal $k$-varieties. 
For a prime divisor $F$ on $X$, 
the {\em push-forward} $f_*F$ of $F$ by $f$ (or to $Y$) is defined as follows: 
if $F \subset \Ex(f)$, then $f_*F=0$, and if $F \not\subset \Ex(f)$, 
then $f_*F$ is defined as the prime divisor whose generic point is equals to the generic point of $f(F)$. 
For an $\R$-divisor $D$ on $X$, 
the {\em push-forward} $f_*D$ of $D$ by $f$ (or to $Y$) 
is defined as $\sum_{i \in I} r_i f_*D_i$, where $D=\sum_{i \in I} r_i D_i$ is the irreducible decomposition of $D$. 
\item 
Let $k$ be a field and 
let $f:X \dashrightarrow Y$ be a birational map. 
Let $D$ be an $\R$-divisor on $X$. 
The {\em proper transform} $D'$ of $D$ on $Y$ is defined as $f'_*D|_{X'}$, 
where $X'$ denotes the maximum open subset $X'$ of $X$ where $f$ is defined, and 
$f':X' \to Y$ is the induced birational morphism. 
\item 
Let $D$ be a closed subset of a smooth scheme $X$ over a perfect field $k$ and 
let $D_1, \ldots , D_n$ be the irreducible components of $D$ with the reduced scheme structures. 
We say that $D$ is \textit{simple normal crossing} if the scheme-theoretic intersection $\bigcap_{j \in J} D_j$  
is smooth over $k$ for every subset $J \subset \{1, \ldots , n \}$. 
For a variety $X$ over $k$ and a closed subset $Z$ of $X$, 
we say that $f:Y \to X$ is a \textit{log resolution} of $(X, Z)$ if $f$ is a projective birational $k$-morphism from a smooth variety $Y$ over $k$ such that 
$\Ex (f) \cup f^{-1}(Z)$ is purely codimension one and 
$\Ex (f) \cup f^{-1}(Z)$ is simple normal crossing. 
In dimension three, there exists a log resolution for such a pair $(X,Z)$ by \cite{CP08}. 
For a variety $X$ over $k$ and an $\R$-divisor $\Delta$ on $X$, 
a log resolution of $(X, \Delta)$ means a log resolution of $(X, \Supp\,\Delta)$. 
\end{itemize}

\subsection{Log pairs}\label{section:log_pair}

We recall some notation in the theory of singularities in the minimal model program. 
For more details, we refer the reader to \cite[Section~2.3]{KM98} and \cite[Section~1]{Kol13}.

We say that $(X, \Delta)$ is a \textit{log pair} over a field $k$ 
if $X$ is a normal variety over $k$ and $\Delta$ is an effective $\R$-divisor such that 
$K_X + \Delta$ is $\R$-Cartier. 
For a proper birational morphism $f: X' \to X$ from a normal variety $X'$ 
and a prime divisor $E$ on $X'$, the \textit{log discrepancy} of $(X, \Delta)$ 
at $E$ is defined as
\[
a_E (X, \Delta) := \coeff _E (K_{X'} - f^* (K_X + \Delta))+1. 
\]

We say that a log pair $(X, \Delta)$ is \textit{log canonical} 
if $a_E (X, \Delta) \ge 0$ for any prime divisor $E$ over $X$. 
Moreover, we say that a log pair $(X, \Delta)$ is \textit{klt}, if $a_E (X, \Delta) > 0$ for any prime divisor $E$ over $X$. 

We call a log pair $(X, \Delta)$ \textit{dlt} 
when all coefficients of $\Delta$ belong to $[0,1]$ and there exists a log resolution $f : Y \to X$ of $(X, \Delta)$ such that 
$a_E(X, \Delta) > 0$ for any $f$-exceptional divisor $E$ on $Y$.

For a log pair $(X, \Delta)$, 
$\Nklt (X, \Delta)$ denotes the closed subset of $X$ 
consisting of the non-klt points of $(X, \Delta)$. 
We equip $\Nklt (X, \Delta)$ with the reduced scheme structure.

\begin{dfn}\label{d-MFS}
Given a field $k$, a log pair $(X, \Delta)$ over $k$, 
and a projective $k$-morphisms $X \xrightarrow{f_1} Z_1 \to Z_2$ 
of quasi-projective $k$-schemes, 
we say that $f_1:X \to Z_1$ is a $(K_X+\Delta)$-{\em Mori fibre space over} $Z_2$ 
if 
$\dim X>\dim Z_1$, $(f_1)_*\MO_X=\MO_{Z_1}$, $\rho(X/Z_1)=1$ and $-(K_X+\Delta)$ is $f_1$-ample. 
If $Z_2=\Spec\,k$, then $f_1:X \to Z_1$ is simply called a $(K_X+\Delta)$-{\em Mori fibre space}. 
\end{dfn}

\subsection{Log minimal models}

\begin{dfn}\label{d-lmm}
Let $k$ be a field and let $Z$ be a quasi-projective $k$-scheme. 
Let $(X, \Delta)$ and $(Y, \Delta_Y)$ be log pairs over $k$ that are projective over $Z$. 
\begin{enumerate}
\item 
We say that $(Y, \Delta_Y)$ is a {\em log birational model over} $Z$ of $(X, \Delta)$ 
if there exists a birational map $\varphi:X \dashrightarrow Y$ over $Z$ 
such that $\Delta_Y=\widetilde \Delta+E$, 
where $\widetilde \Delta$ denotes the proper transform of $\Delta$ on $Y$ and 
$E$ is the sum of the prime divisors that are exceptional over $X$. 
\item 
We say that $(Y, \Delta_Y)$ is a {\em weak log canonical model over} $Z$ of $(X, \Delta)$ 
if 
\begin{enumerate}
\item 
$(Y, \Delta_Y)$ is a log birational model over $Z$ of $(X, \Delta)$, 
\item 
$K_Y+\Delta_Y$ is nef over $Z$, and 
\item 
for any prime divisor $D$ on $X$ that is exceptional over $Y$, it holds that 
\[
a_D(X, \Delta) \leq a_D(Y, \Delta_Y). 
\]
\end{enumerate}
\item 
We say that $(Y, \Delta_Y)$ is a {\em log canonical model over} $Z$ of $(X, \Delta)$ 
if 
\begin{enumerate}
\item 
$(Y, \Delta_Y)$ is a weak log canonical model over $Z$ of $(X, \Delta)$, and
\item 
$K_Y+\Delta_Y$ is ample over $Z$. 
\end{enumerate}
\item 
We say that $(Y, \Delta_Y)$ is a {\em log minimal model over} $Z$ of $(X, \Delta)$ 
if 
\begin{enumerate}
\item[(i)] 
$(Y, \Delta_Y)$ is a weak log canonical model over $Z$ of $(X, \Delta)$, 
\item[(ii)]  
$(Y, \Delta_Y)$ is a $\Q$-factorial dlt pair, and 
\item[(iii)]  
for any prime divisor $D$ on $X$ that is exceptional over $Y$, it holds that 
\[
a_D(X, \Delta) < a_D(Y, \Delta_Y). 
\]
\end{enumerate}
\end{enumerate}
\end{dfn}

\begin{rem}\label{r-2-wlc-models}
Let $k$ be a field and let $Z$ be a quasi-projective $k$-scheme. 
Fix a log canonical pair $(X, \Delta)$ over $k$ that is projective over $Z$. 
Let $(Y_1, \Delta_{Y_1})$ and $(Y_2, \Delta_{Y_2})$ be weak log canonical models over $Z$ 
of $(X, \Delta)$. 
If $g_1:W \to Y_1$ and $g_2:W \to Y_2$ are projective birational morphisms 
that commute with $Y_1 \dashrightarrow Y_2$, 
then it follows from the same argument as in \cite[Remark 2.6]{Bir12} that 
\[
g_1^*(K_{Y_1}+\Delta_{Y_1})=g_2^*(K_{Y_2}+\Delta_{Y_2}). 
\]
Moreover, if $(Y_2, \Delta_{Y_2})$ is a log canonical model of $(X, \Delta)$, 
then $K_{Y_1}+\Delta_{Y_1}$ is semi-ample over $Z$ and 
the birational map $Y_1 \dashrightarrow Y_2$ is a morphism. 
In particular, if both $(Y_1, \Delta_{Y_1})$ and $(Y_2, \Delta_{Y_2})$ are 
log canonical models of $(X, \Delta)$, then the induced birational map $Y_1 \dashrightarrow Y_2$ 
is an isomorphism.  
\end{rem}

\subsection{Minimal model program}

%\subsubsection{\textcolor{red}{Definition}}

Let us recall construction of a sequence of steps of log canonical minimal model program.

\begin{dfn}[{cf.\ \cite[4.9.1]{Fuj17}}]
Let $k$ be a field. 
Let $(X, \Delta)$ be a log pair over $k$ and let $f:X \to Z$ be 
a projective $k$-morphism to a quasi-projective $k$-scheme $Z$. 
Let $g:X \dashrightarrow Y$ be a birational map over $Z$ to a normal $k$-variety $Y$ projective over $Z$. 
\begin{enumerate}
\item 
We say that $g$ is a $(K_X+\Delta)$-{\em divisorial contraction} over $Z$ if 
\begin{enumerate}
\item $g$ is a morphism, 
\item $\dim \Ex(g)=\dim X-1$, 
\item $-(K_X+\Delta)$ is $g$-ample, and 
\item $\rho(X/Y)=1$. 
\end{enumerate}
\item 
We say that $g$ is a $(K_X+\Delta)$-{\em flipping contraction} over $Z$ if 
\begin{enumerate}
\item $g$ is a morphism, 
\item $\dim \Ex(g)<\dim X-1$, 
\item $-(K_X+\Delta)$ is $g$-ample, and 
\item $\rho(X/Y)=1$. 
\end{enumerate}
\item 
We say that $g$ is a {\em step of a }$(K_X+\Delta)$-{\em MMP} over $Z$ if 
there exist a $(K_X+\Delta)$-divisorial or $(K_X+\Delta)$-flipping contraction $h:X \to W$ over $Z$ 
and a log canonical model $(Y=X^+, \Delta^+)$ of $(X, \Delta)$ over $W$ such that 
$g=(h^+)^{-1} \circ h$, where 
$h^+:X^+ \to W$ is the induced morphism. 
If $h$ is a flipping contraction, then any of $h^+$ and $g$ is called a $(K_X+\Delta)$-{\em flip} over $Z$. 
\end{enumerate}
\end{dfn}

\begin{dfn}
Let $k$ be a field. 
Let $(X, \Delta)$ be a log pair over $k$ and let $f:X \to Z$ be 
a projective $k$-morphism to a quasi-projective $k$-scheme $Z$. 
\begin{enumerate}
\item 
A (possibly infinite) sequence 
\[
X=X_0 \overset{\varphi_0}{\dashrightarrow} X_1 \overset{\varphi_1}{\dashrightarrow} \cdots 
\]
is called a $(K_X+\Delta)$-{\em MMP over}  $Z$ if 
any $\varphi_i$ is a step of a $(K_X+\Delta)$-MMP over $Z$. 
If $Z=\Spec\,k$, then a $(K_X+\Delta)$-MMP over  $Z$ 
is simply called a $(K_X+\Delta)$\textit{-MMP}.  
\item 
For an $\R$-Cartier $\R$-divisor $C$ on $X$, a (possibly infinite) sequence 
\[
X=X_0 \overset{\varphi_0}{\dashrightarrow} X_1 \overset{\varphi_1}{\dashrightarrow} \cdots 
\]
is called a $(K_X+\Delta)$-{\em MMP over} $Z$ {\em with scaling of} $C$ if 
the sequence is a $(K_X+\Delta)$-MMP over $Z$ and 
\[
\lambda_i :=\min\{\lambda \in \R_{\geq 0} \,|\, K_{X_i}+\Delta_i +\lambda C_i\text{ is nef over }Z\}
\]
exists and $K_{X_i}+\Delta_i +\lambda_i C_i$ is numerically trivial over $Z_i$, 
where $X_i \to Z_i$ is the associated divisorial or flipping contraction. 
In this case, $\lambda_0, \lambda_1, \cdots$ are called the {\em scaling coefficients}. 
Note that we do not assume that $C$ is effective. 
\item 
We say that a $(K_X+\Delta)$-MMP over $Z$ {\em terminates} 
if the sequence is a finite sequence: 
\[
X=X_0 \overset{\varphi_0}{\dashrightarrow} X_1 \overset{\varphi_1}{\dashrightarrow} 
\cdots \overset{\varphi_{\ell-1}}{\dashrightarrow} X_{\ell}
\]
such that either $K_{X_{\ell}}+\Delta_{\ell}$ is nef over $Z$ or 
there is a $(K_{X_{\ell}}+\Delta_{\ell})$-Mori fibre space over $Z$. 
\end{enumerate}
\end{dfn}

%\subsubsection{\textcolor{red}{Known results on log canonical minimal model program}}

For later use, we now summarise known results from \cite{Wal} on three-dimensional log canonical minimal model program 
in positive characteristic.

\begin{thm}\label{t-Wal1.1}
Let $k$ be an algebraically closed field of characteristic $p>5$. 
Let $(X, \Delta)$ be a three-dimensional log canonical pair over $k$ 
and let $f:X \to Z$ be a projective $k$-morphism to a projective $k$-scheme $Z$. 
If $K_X+\Delta$ is $f$-nef and $f$-big, then $K_X+\Delta$ is $f$-semi-ample. 

\end{thm}

\begin{proof}
See \cite[Theorem 1.1]{Wal}. 
\end{proof}

\begin{thm}\label{t-Wal1.2}
Let $k$ be an algebraically closed field of characteristic $p>5$. 
Let $(X, \Delta)$ be a three-dimensional log canonical pair over $k$ 
and let $f:X \to Z$ be a projective $k$-morphism to a projective $k$-scheme $Z$. 
If $A$ is an $f$-semi-ample and $f$-big $\R$-Cartier $\R$-divisor 
such that $K_X+\Delta+A$ is $f$-nef, 
then $K_X+\Delta+A$ is $f$-semi-ample. 
\end{thm}

\begin{proof}
See \cite[Theorem 1.2]{Wal}. 
\end{proof}

\begin{thm}\label{t-Wal1.7}
Let $k$ be an algebraically closed field of characteristic $p>5$. 
Let $(X, \Delta)$ be a three-dimensional projective log canonical pair over $k$. 
Then there exists a countable set $\{ C_i \} _{i \in I}$ 
of rational curves on $X$ which satisfies the following conditions: 
\begin{enumerate}
\item $\overline{\NE}(X)=\overline{\NE}(X)_{K_X+\Delta \geq 0}+\sum_{i \in I} \R_{\geq 0}  [C_i]$. 
\item $0 < - (K_X + \Delta) \cdot C_i \le 6$ for any $i \in I$. 
\item 
For any ample $\R$-Cartier $\mathbb{R}$-divisor $A$, 
there exists a finite subset $J$ of $I$ such that 
\[
\overline{\NE}(X)=\overline{\NE}(X)_{K_X+\Delta+A \geq 0}+\sum_{j \in J} \R_{\geq 0}  [C_j].
\] 
\end{enumerate}
\end{thm}

\begin{proof}
See \cite[Theorem 1.7]{Wal}. 
\end{proof}

\begin{thm}\label{t-Wal1.8}
Let $k$ be an algebraically closed field of characteristic $p>5$. 
Let $(X, \Delta)$ be a three-dimensional $\Q$-factorial log canonical pair over $k$ 
and let $f:X \to Z$ be a projective $k$-morphism to a projective $k$-scheme $Z$. 
Assume that there exists an effective $\R$-divisor $M$ on $X$ 
such that $K_X + \Delta \equiv_f M$. 
Then there exists a $(K_X+\Delta)$-MMP over $Z$ that terminates. 
\end{thm}

\begin{proof}
See \cite[Corollary 1.8]{Wal}. 
\end{proof}

\subsection{Compactifications of pairs}

First we show the existence of a compactification of klt pairs. 

\begin{prop}\label{p-klt-cpt}
Let $k$ be a perfect field of characteristic $p>5$. 
Let $(X, \Delta)$ be a three-dimensional quasi-projective $\Q$-factorial klt pair over $k$ and let $f:X \to Z$ be a $k$-morphism to a quasi-projective $k$-scheme $Z$. 
Then there exists an open immersion $j:X \to \overline X$ over $Z$ 
to a normal $\Q$-factorial threefold $\overline X$ projective over $Z$ such that 
$(\overline X, \overline{\Delta})$ is klt for the closure $\overline{\Delta}$ of $\Delta$. 
\end{prop}

\begin{proof}
Enlarging coefficients of $\Delta$ appropriately, 
we may assume that $\Delta$ is a $\Q$-divisor. 
Since $X$ is quasi-projective over $Z$, 
there exists an open immersion $j_1:X \to X_1$ over $Z$ 
to an integral normal scheme $X_1$ projective over $Z$. 

Let $\varphi:V \to X_1$ be a resolution of singularities of $X_1$ 
such that $\varphi^{-1}(\overline{\Supp\,\Delta}) \cup \Ex(\varphi)$ is 
a simple normal crossing divisor. 
Let $E$ be the reduced divisor on $V$ such that $\Supp\,E=\Ex(\varphi)$. 
Fix sufficiently small $\epsilon \in \Q_{>0}$. 
Then it holds that $(V, \varphi_*^{-1}\overline{\Delta}+(1-\epsilon)E)$ is klt. 
Moreover, it follows from \cite[Theorem 2.13]{GNT} 
that there exists a $(K_V+\varphi_*^{-1}\overline{\Delta}+(1-\epsilon)E)$-MMP over $X_1$ 
that terminates. 
Let $\overline X$ be the end result. 
Then the negativity lemma and the $\mathbb{Q}$-factoriality of $X$ imply that 
$\overline X \to X_1$ is isomorphic over $X$. 
Moreover, $(\overline X, \overline{\Delta})$ is klt 
for the closure $\overline{\Delta}$ of $\Delta$ in $\overline{X}$. 
\end{proof}

Using Proposition \ref{p-klt-cpt} and the ACC for log canonical thresholds, 
we prove the existence of a certain compactification of log canonical pairs 
(Proposition \ref{p-lc-cpt}, Proposition \ref{p-lc-cpt3}).

\begin{prop}\label{p-lc-cpt}
Let $k$ be a perfect field of characteristic $p>5$. 
Let $(X, \Delta)$ be a three-dimensional quasi-projective $\Q$-factorial log canonical pair over $k$ such that $(X, \{\Delta\})$ is klt. 
Let $f:X \to Z$ be a $k$-morphism to a quasi-projective $k$-scheme $Z$. 
Then there exists an open immersion $j:X \to \overline X$ over $Z$
to a normal $\Q$-factorial threefold $\overline X$ projective over $Z$ 
such that $(\overline X, \overline \Delta)$ is log canonical and 
$(\overline X, \{\overline \Delta\})$ is klt, 
where $\overline \Delta$ denotes the closure of $\Delta$ in $\overline X$. 
\end{prop}

\begin{proof}
Since $(X, \{\Delta\})$ is klt, also $(X, \Delta-\epsilon \llcorner \Delta\lrcorner)$ 
is klt for any $\epsilon \in \R$ such that $0<\epsilon\leq 1$. 
For a sufficiently small $\epsilon \in \R_{>0}$, 
we apply Proposition \ref{p-klt-cpt} to a $\Q$-factorial klt pair 
$(X, \Delta-\epsilon \llcorner \Delta\lrcorner)$ and a morphism $X \to Z$. 
Then we get an open immersion $j:X \to \overline X$ over $Z$ 
to a normal $\Q$-factorial threefold $\overline X$ projective over $Z$ such that 
$(\overline X, D)$ is klt for the closure $D$ of $\Delta-\epsilon \llcorner \Delta\lrcorner$. 
Since $\epsilon \in \R_{>0}$ is sufficiently small, 
it follows from the ACC for log canonical thresholds \cite[Theorem 1.10]{Bir16}
that $(\overline X, \overline \Delta)$ is log canonical for the closure $\overline \Delta$ 
of $\Delta$. 
\end{proof}

\begin{prop}\label{p-lc-cpt3}
Let $k$ be a perfect field of characteristic $p>5$. 
Let $(X, \Delta)$ be a three-dimensional quasi-projective $\Q$-factorial log canonical pair over $k$ such that $X$ is klt. 
Let $f:X \to Z$ be a $k$-morphism to a quasi-projective $k$-scheme $Z$. 
Then there exists an open immersion $j:X \to \overline X$ over $Z$
to a normal $\Q$-factorial threefold $\overline X$ projective over $Z$ 
such that $(\overline X, \overline \Delta)$ is log canonical and 
$\overline X$ is klt. 
\end{prop}

\begin{proof}
Since $(X, 0)$ is klt, also $(X, (1-\epsilon)\Delta)$ 
is klt for any $\epsilon \in \R$ such that $0<\epsilon\leq 1$. 
For a sufficiently small $\epsilon \in \R_{>0}$, 
we apply Proposition \ref{p-klt-cpt} to a $\Q$-factorial klt pair 
$(X, (1-\epsilon)\Delta)$ and a morphism $X \to Z$. 
Then we get an open immersion $j:X \to \overline X$ over $Z$ 
to a normal $\Q$-factorial threefold $\overline X$ projective over $Z$ such that 
$(\overline X, D)$ is klt for the closure $D$ of $(1-\epsilon)\Delta$. 
Since $\epsilon \in \R_{>0}$ is sufficiently small, 
it follows from the ACC for log canonical thresholds \cite[Theorem 1.10]{Bir16}
that $(\overline X, \overline \Delta)$ is log canonical for the closure $\overline \Delta$ 
of $\Delta$. 
\end{proof}

\subsection{Perturbation of coefficients}

\begin{lem}\label{l-lcklt-perturb}
Fix $\mathbb K \in \{\Q, \R\}$. 
Let $k$ be an infinite perfect field of characteristic $p>5$. 
Let $(X, \Delta)$ be a three-dimensional log pair over $k$ and 
let $f:X \to Z$ be a projective $k$-morphism to a quasi-projective $k$-scheme. 
Let $A$ be an $f$-semi-ample $\mathbb K$-Cartier $\mathbb K$-divisor on $X$. 
Then the following hold. 
\begin{itemize}
\item[(1)] If $(X, \Delta)$ is log canonical, then there exists an effective $\mathbb K$-Cartier 
$\mathbb K$-divisor $A'$ such that 
$A \sim_{f, \mathbb K} A'$ and $(X, \Delta+A')$ is log canonical. 
\item[(2)] If $(X, \Delta)$ is klt, then there exists an effective $\mathbb K$-Cartier 
$\mathbb K$-divisor $A'$ such that 
$A \sim_{f, \mathbb K} A'$ and $(X, \Delta+A')$ is klt. 
\end{itemize}
\end{lem}

\begin{proof}
First we prove (1).
By standard argument, we may assume that $A$ is an $f$-semi-ample $\mathbb Q$-Cartier $\mathbb Q$-divisor (note that $\Delta$ could be still an $\R$-divisor). 
Taking a log resolution of $(X, \Delta)$, we may assume that $X$ is smooth over $k$ and 
$\Delta$ is a simple normal crossing $\R$-divisor. 
Hence, enlarging coefficients of $\Delta$, the problem is reduced to the case when $\mathbb K=\Q$. 

Since $A$ is an $f$-semi-ample $\Q$-divisor, there exist projective morphisms of schemes 
\[
f:X \xrightarrow{g} Y \xrightarrow{h} Z
\]
and an $h$-ample $\Q$-Cartier $\Q$-divisor $A_Y$ 
such that $Y$ is a normal variety over $k$ and $A \sim_{\Q} g^*A_Y$. 

Fix an open immersion $j_Z:Z \hookrightarrow \overline Z$ 
to a projective $k$-scheme. 
We can find an open immersion $j_Y:Y \hookrightarrow \overline Y$ over $\overline Z$ 
to a normal $k$-variety $\overline Y$ and 
a globally ample $\Q$-Cartier $\Q$-divisor $A_{\overline Y}$ such that 
$\overline Y$ is projective over $\overline Z$ and 
$A_{\overline Y}|_Y \sim_{Z, \Q} A_Y$.  
Applying Proposition \ref{p-lc-cpt3} to a log canonical pair $(X, \Delta)$ and a morphism 
$X \to \overline Y$, 
there exists an open immersion $j_X:X \hookrightarrow \overline X$ over $\overline Z$ 
such that $\overline X$ is a $\Q$-factorial threefold over $k$ and 
$(\overline X, \overline \Delta)$ is log canonical, where $\overline{\Delta}$ denotes the closure of 
$\Delta$ in $\overline X$. 
To summarise, we have a commutative diagram: 
$$\begin{CD}
X @>j_X >> \overline X\\
@VVgV @VV\overline g V\\
Y @>j_Y >> \overline Y\\
@VVhV @VV\overline h V\\
Z @>j_Z >> \overline Z. 
\end{CD}$$

We apply \cite[Theorem 1]{Tanb} to a log canonical pair $(\overline X, \overline \Delta)$ and a semi-ample $\Q$-divisor $\overline g^*A_{\overline Y}$. 
Then there exists an effective $\Q$-divisor $B$ on $\overline X$ such that 
$\overline{g}^*A_{\overline Y} \sim_{\Q} B$ and $(\overline X, \overline \Delta+B)$ is log canonical. 
Set $A':=B|_X$. 
It holds that $(X, \Delta+A')$ is log canonical and 
\[
A'=B|_X \sim_{\Q}  (\overline{g}^*A_{\overline Y})|_X \sim_{Z, \Q} g^*A_Y \sim_{Z, \Q} A.
\]
We have proved (1).

Next we prove (2). Suppose that $(X, \Delta)$ is klt. 
We apply (1) to $(X, \Delta)$ and $2A$. 
Then there exists an effective $\mathbb K$-Cartier $\mathbb K$-divisor $B$ such that 
$2A \sim_{f, \mathbb K} B$ and $(X, \Delta+B)$ is log canonical. 
Set $A':=(1/2)B$. 
Then we have that $A \sim_{f, \mathbb K} A'$ and $(X, \Delta+A')$ is klt. 
\end{proof}

\begin{lem}\label{l-coeff-1}
Let $k$ be a perfect field. 
Let $X$ be a quasi-projective $\Q$-factorial normal variety over $k$. 
Let $\Delta$ be an effective $\Q$-divisor such that all the coefficients of $\Delta$ are contained in $[0, 1]$. 
Let $A$ be an ample $\Q$-divisor on $X$. 
Then there exists an effective $\Q$-divisor $\Delta'$ on $X$ such that 
$\Delta' \sim_{\Q} \Delta+A$ and all the coefficients of $\Delta'$ are contained in $[0, 1]$. 
\end{lem}

\begin{proof}
There exists a positive real number $\epsilon$ such that $0< \epsilon < 1$ and 
$\epsilon \Delta + A$ is ample. 
After replacing $\Delta$ and $A$ by $(1-\epsilon)\Delta$ and $\epsilon A$ respectively, 
we may assume that there exists a positive real number $\epsilon$ such that $0< \epsilon < 1$ and 
all the coefficients of $\Delta$ are contained in $[0, 1-\epsilon]$. 
%Fix a positive integer $n_1$ such that $1/n_1<\epsilon$. 
If $k$ is an infinite field (resp.\ a finite field), 
then \cite[Theorem 1]{Sei50} (resp.\ \cite[Theorem 1.1]{Poo04}) enables us to find 
a positive integer $n$ and a smooth divisor $H$ on 
the smooth locus $X_{{\rm sm}}$ of $X$ such that 
\begin{itemize}
\item $1/n<\epsilon$, and 
\item $nA$ is a Cartier divisor such that $(n A)|_{X_{{\rm sm}}} \sim H$.
\end{itemize}
For the closure $\overline H$ of $H$ in $X$, it holds that 
$nA \sim \overline H$ and $\overline H$ is a reduced divisor. 
Then it holds that 
\[
\Delta+A=\Delta+\frac{1}{n}(nA) \sim_{\Q} \Delta+\frac{1}{n} \overline H=:\Delta'
\]
and all the coefficients of $\Delta'$ are contained in $[0, 1]$. 
\end{proof}

\subsection{Non-vanishing theorem of relative dimension two}

\begin{lem}\label{l-2dim-nonvani}
Let $k$ be a field of characteristic $p>0$. 
Let $(X, \Delta)$ be a log canonical pair over $k$ and 
let $f:X \to Z$ be a projective $k$-morphism to a quasi-projective $k$-variety such that $f_*\MO_X=\MO_Z$. 
Set $X_K$ to be the generic fibre of $f$. 
Assume that 
\begin{enumerate}
\item $\dim X_K \leq 2$, and 
\item $(K_X+\Delta)|_{X_K}$ is pseudo-effective. 
\end{enumerate}
Then there exists an effective $\R$-Cartier $\R$-divisor $D$ on $X$ 
such that $K_X+\Delta \sim_{Z, \R} D$. 
\end{lem}

\begin{proof}
By the non-vanishing theorem for surfaces (cf.\ \cite[Theorem 1.1]{Tanc} and \cite[Theorem 1.1]{Tand}), 
we obtain 
\[
(K_X+\Delta)|_{X_K}+\sum_{i=1}^a r_i {\rm div}(\varphi_i)=E
\]
for some effective $\R$-divisor $E$ on $X_K$, 
$r_1, \ldots, r_a \in \R$, and $\varphi_1, \ldots, \varphi_a \in K(X_K)$. 
By $K(X)=K(X_K)$, we get $\varphi_i \in K(X)$. 
We define an $\R$-Cartier $\R$-divisor $F$ by 
\[
K_X+\Delta+\sum_{i=1}^a r_i {\rm div}(\varphi_i)=:F.
\]
We obtain $F|_{X_K}=E.$ 
In particular, if $F=F_1-F_2$ for effective $\R$-divisors $F_1$ and $F_2$ 
without any common irreducible components, then 
we get $f(\Supp\,F_2) \subsetneq Z$. 
Since $Z$ is quasi-projective over $k$, 
there exists an effective Cartier divisor $H_Z$ on $Z$ such that 
$f(\Supp\,F_2) \subset \Supp\,H_Z$. 
Hence, for a sufficiently large positive integer $m$, it holds that 
$mf^*H_Z-F_2$ is an effective $\R$-divisor. 
We have that 
\[
K_X+\Delta+\sum_{i=1}^a r_i {\rm div}(\varphi_i)=F=F_1-F_2 
\sim_{Z, \R} F_1+(mf^*H_Z-F_2) \ge 0,
\]
as desired. 
\end{proof}

\section{MMP for the effective case} 

The purpose of this section is to establish the relative minimal model program 
for log canonical threefolds with effective log canonical divisors (Theorem \ref{t-lc-eff-mmp}). 
As a consequence, 
we obtain the existence of dlt modifications (Proposition \ref{p-dlt-modif}). 
We start with the following lemma on dlt modification. 
Note that this result is known if the base field is algebraically closed \cite[Theorem 1.6]{Bir16}.

\begin{lem}\label{l-lc-modif}
Let $k$ be a perfect field of characteristic $p>5$. 
Let $(X, \Delta)$ be a three-dimensional quasi-projective $\Q$-factorial log canonical pair over $k$, 
where $\Delta$ is a $\Q$-divisor. 
Then there exists a projective birational $k$-morphism $f:Y \to X$ 
from a normal threefold $Y$ over $k$ such that 
\begin{enumerate}
\item $Y$ is $\Q$-factorial,
\item $a_F(X, \Delta)=0$ for any $f$-exceptional prime divisor $F$ on $Y$, and 
\item if $\Delta_Y$ is the $\R$-divisor defined by $K_Y+\Delta_Y=f^*(K_X+\Delta)$, then 
$\Delta_Y$ is effective and $(Y, \{\Delta_Y\})$ is klt 
(see Subsection \ref{ss-notation} for the definition of $\{\Delta_Y\}$). 
\end{enumerate}
\end{lem}

\begin{proof}
Let $\varphi:V \to X$ be a log resolution of $(X, \Delta)$. 
We can write 
$$K_V+\varphi^{-1}_*\Delta+E=\varphi^*(K_X+\Delta)+\sum_i a_i E_i,$$
where $E_i$ is a $\varphi$-exceptional prime divisor, 
$a_i:=a_{E_i}(X, \Delta) \in \Q_{\geq 0}$ is the log discrepancy, and 
$E$ is the reduced divisor on $V$ such that $\Supp\,E=\Ex(\varphi)$. 
Fix a sufficiently small $\epsilon \in \Q_{>0}$ so that 
if $a_i > 0$, then $a_i-\epsilon>0$. 
Since $\varphi^{-1}_*\Delta+(1-\epsilon)E - \delta \varphi ^{*} \Delta$
is still effective for sufficiently small $\delta >0$, 
\cite[Theorem 2.13]{GNT} enables us to find  
a $(K_V+\varphi^{-1}_*\Delta+(1-\epsilon)E)$-MMP over $X$ that terminates. 
Let $g:Y \to X$ be the end result. 
Then (1) holds. 
The negativity lemma implies that this MMP contracts 
all the prime divisors $E_i$ satisfying $a_i=a_{E_i}(X, \Delta)>0$. 
In other words, any $g$-exceptional prime divisor $F$ satisfies $a_F(X, \Delta)=0$. 
Therefore, (2) holds. 
Furthermore, (3) holds because being dlt is preserved under the MMP. 
\end{proof}

\begin{thm}\label{t-lc-flip}
Let $k$ be a perfect field of characteristic $p>5$. 
Let $(X, \Delta)$ be a three-dimensional $\Q$-factorial log canonical pair over $k$, where $\Delta$ is a $\Q$-divisor. Let $f:X \to Z$ be a projective $k$-morphism to a quasi-projective $k$-scheme $Z$. 
If $K_X+\Delta$ is $f$-big, then the graded ring
$$\bigoplus_{m=0}^{\infty} f_*\MO_X(\llcorner m(K_X+\Delta)\lrcorner)$$
is a finitely generated $\MO_Z$-algebra. 
\end{thm}

\begin{proof}
By standard argument, 
we may assume that $k$ is algebraically closed and $f_*\MO_X=\MO_Z$ (cf.\ \cite[Remark 2.7]{GNT}). 
By Lemma \ref{l-lc-modif}, the problem is reduced to the case when 
\begin{enumerate}
\item 
$(X, \{\Delta\})$ is klt. 
\end{enumerate}
Take an open immersion $j_Z:Z \to \overline Z$ 
such that $\overline Z$ is a projective normal threefold over $k$. 
Applying Proposition \ref{p-lc-cpt} to $(X, \Delta)$ and $X \to \overline Z$, 
we may assume that (1) and 
\begin{enumerate}
\item[(2)]  
both $X$ and $Z$ are projective over $k$. 
\end{enumerate}
By Theorem \ref{t-Wal1.8}, there exists a $(K_X+\Delta)$-MMP over $Z$ that terminates. 
Replacing $(X, \Delta)$ by the end result, 
we may assume that (2) and 
\begin{enumerate}
\item[(3)]  
$K_X+\Delta$ is $f$-nef. 
\end{enumerate}
By Theorem \ref{t-Wal1.2}, %\redout{\cite[Theorem 1.2]{Wal}}, $K_X+\Delta$ is $f$-semi-ample. 
%Therefore, 
the assertion in the statement holds. 
\end{proof}

\begin{prop}\label{p-lc-big-bpf}
Let $k$ be a perfect field of characteristic $p>5$. 
Let $(X, \Delta)$ be a three-dimensional $\Q$-factorial log canonical pair over $k$, where $\Delta$ is a $\Q$-divisor. Let $f:X \to Z$ be a projective $k$-morphism to a quasi-projective $k$-scheme $Z$. 
Let $L$ be an $f$-nef and $f$-big $\Q$-Cartier $\Q$-divisor such that 
$L-(K_X+\Delta)$ is $f$-semi-ample. 
Then $L$ is $f$-semi-ample. 
\end{prop}

\begin{proof}
We may assume that $k$ is algebraically closed and $f_*\MO_X=\MO_Z$. 
By Lemma \ref{l-lcklt-perturb} (1), 
we are reduced to the case when $L = K_X+\Delta$. 
By Lemma \ref{l-lc-modif}, 
we may assume that $(X, \{\Delta\})$ is klt. 

Fix an open immersion $j_Z:Z \hookrightarrow \overline Z$ 
such that $\overline Z$ is a normal variety projective over $k$. 
It follows from Proposition \ref{p-lc-cpt} that 
there exists an open immersion $j_1:X \hookrightarrow X_1$ over $\overline Z$ 
to a $\Q$-factorial threefold $X_1$ projective over $\overline Z$ such that 
$(X_1, \Delta_{X_1})$ is log canonical for the closure $\Delta_{X_1}$ of $\Delta$ in $X_1$. 
Let $f_1:X_1 \to \overline Z$ be the induced morphism. 
Since $K_{X_1}+\Delta_{X_1} \sim_{f_1, \Q} D$ for some effective $\Q$-divisor $D$, 
Theorem \ref{t-Wal1.8} enables us to find 
a $(K_{X_1}+\Delta_{X_1})$-MMP over $\overline Z$ that terminates. 
Let $X_2$ be the end result and let $\Delta_{X_2}$ be the push-forward of $\Delta_{X_1}$. 
Then $K_{X_2}+\Delta_{X_2}$ is $f_2$-nef, 
where $f_2:X_2 \to \overline Z$ is the induced morphism. 
Since $K_X+\Delta$ is $f$-nef, the two morphisms
$f_1$ and $f_2$ coincide over $Z$. 
In particular, there exists an open immersion $j_2:X \hookrightarrow X_2$ 
such that $j_2(X)=f_2^{-1}(Z)$. 
It follows from Theorem \ref{t-Wal1.2} that $K_{X_2}+\Delta_{X_2}$ is $f_2$-semi-ample. 
Taking the restriction to $X$, we see that $K_X+\Delta$ is $f$-semi-ample, as desired. 
\end{proof}

\begin{thm}\label{t-lc-eff-mmp}
Let $k$ be a perfect field of characteristic $p>5$. 
Let $(X, \Delta)$ be a three-dimensional $\Q$-factorial log canonical pair over $k$ and 
let $f:X \to Z$ be a projective $k$-morphism to a quasi-projective $k$-scheme $Z$. 
Assume that $K_X+\Delta \equiv_{f, \R} D$ for some effective $\R$-divisor $D$ on $X$. 
Then there exists a $(K_X+\Delta)$-MMP over $Z$ that terminates. 
\end{thm}

\begin{proof}
Assuming that $K_X+\Delta$ is not $f$-nef, 
let us find a $(K_X+\Delta)$-negative extremal ray and its contraction. 
There exists an $f$-ample $\R$-divisor $A$ such that 
$K_X+\Delta+A$ is not $f$-nef.

Let us find an effective $\Q$-divisor $\Delta'$, 
an $f$-ample $\Q$-divisor $A'_1$ and an $f$-ample $\R$-divisor $A'_2$ such that 
\begin{enumerate}
\item $K_X+\Delta+A = K_X+\Delta'+A'_1+A'_2$, 
\item all the coefficients of $\Delta'$ are contained in $[0, 1]$, and 
\item $K_X+\Delta' \equiv_{f, \Q} D'$ for some effective $\Q$-divisor $D'$.  
\end{enumerate}
We first take an ample $\Q$-divisor $A'_1$ on $X$ such that $A-3A'_1$ is $f$-ample. 
Then there exists a $\Q$-divisor $\Delta_1$ on $X$ such that 
$0 \leq \Delta_1 \leq \Delta$ and both $A_1'+(\Delta-\Delta_1)$ 
and $A_1'-(\Delta-\Delta_1)$ are $f$-ample. 
By Lemma \ref{l-coeff-1}, 
there exists an effective $\Q$-divisor $\Delta'$ such that $\Delta' \sim_{\Q} \Delta_1+A'_1$ and all the coefficients of $\Delta'$ are contained in $[0, 1]$. 
Thus (2) holds. 
Then we obtain (3), since the $\R$-divisors appearing in the following:  
\[
(K_X+\Delta')-(K_X+\Delta) \sim_{f, \R} \Delta_1+A_1'-\Delta=A_1'-(\Delta-\Delta_1)
\]
are $f$-ample. 
Set $A'_2:=(K_X+\Delta+A)-(K_X+\Delta'+A'_1)$. 
By the equations: 
\begin{eqnarray*}
A'_2&=&(K_X+\Delta+A)-(K_X+\Delta'+A'_1)\\
&=&(A-3A'_1)+(A'_1+(\Delta-\Delta_1))+(A'_1+(\Delta_1-\Delta')), 
\end{eqnarray*}
it holds that $A'_2$ is $f$-ample, hence (1) holds. 

Applying \cite[Lemma 2.2]{GNT} to $(X, \Delta')$ and $A'_1$, 
there exist finitely many curves $C_1, \ldots, C_r$ on $X$ such that 
\[
\overline{\NE}(X/Z)=\overline{\NE}(X/Z)_{K_X+\Delta '+A'_1+A'_2 \ge 0}+\sum_{i=1}^r \R_{\geq 0}[C_i]. 
\]
Set $R:=\R_{\geq 0}[C_1]$, which is a $(K_X+\Delta)$-negative extremal ray. 
By Proposition \ref{p-lc-big-bpf}, there exists a contraction $\varphi:X \to Y$ of $R$. 
If $\varphi$ is a flipping contraction, then Theorem \ref{t-lc-flip} 
implies that a flip of $\varphi$ exists. 

Therefore, it suffices to prove that there exists no infinite sequence 
that is a $(K_X+\Delta)$-MMP consisting of flips: 
\[
X=X_0 \dashrightarrow X_1 \dashrightarrow \cdots.
\]
This follows from the ACC for log canonical thresholds \cite[Theorem 1.10]{Bir16} and the assumption that 
$K_X+\Delta \equiv_{f, \mathbb{R}} D \ge 0$ (\cite[Lemma 3.2]{Bir07}). 
\end{proof}

Now we prove the existence of dlt modifications (Corollary \ref{c-dlt-modif}). 
First we treat more general pairs in Proposition \ref{p-dlt-modif}. 

\begin{prop}\label{p-dlt-modif} 
Let $k$ be a perfect field of characteristic $p>5$. 
Let $(X, \Delta)$ be a three-dimensional quasi-projective log pair over $k$. 
Then there exists a projective birational $k$-morphism $f:Y \to X$ 
from a normal threefold over $k$ that satisfies the following conditions: 
\begin{enumerate}
\item $a_F (X, \Delta) \le 0$ holds for any $f$-exceptional prime divisor $F$ on $Y$. 
\item $(Y, \Delta _Y ^{\wedge 1})$ is a $\Q$-factorial dlt pair,  
where $\Delta _Y$ is the $\R$-divisor defined  by $K_Y + \Delta _Y = f^* (K_X + \Delta)$ 
(see Subsection \ref{ss-notation} for the definition of $\Delta _Y ^{\wedge 1}$). 
\item 
$\Nklt (Y, \Delta _{Y}) = f ^{-1} (\Nklt (X, \Delta))$ holds. 
\end{enumerate} 
\end{prop}

\begin{proof}
The proof consists of two steps. 

\setcounter{step}{0}
\begin{step}\label{s1-dlt-modif} 
There exists a projective birational morphism $f_1:Y_1 \to X$ that satisfies 
the conditions 
(1) and (2) in the statement. 
\end{step}

\begin{proof}[Proof of Step \ref{s1-dlt-modif}] 
Let $g: W \to X$ be a log resolution of $(X, \Delta)$. Write
\[
g^*(K_X + \Delta) = K_W + \Delta_W = K_W + \widetilde{\Delta} + E, 
\]
where $\widetilde{\Delta}$ is the proper transform of $\Delta$, and 
$E:=\Delta_W-\widetilde{\Delta}$. 
Hence, $E$ is a $g$-exceptional $\mathbb{R}$-divisor. 
Set $F$ to be the reduced divisor such that $\Supp\,F=\Ex(g)$. 
There exists a $(K_W + \widetilde{\Delta} ^{\wedge 1} + F)$-MMP over $X$ 
that terminates (Theorem \ref{t-lc-eff-mmp}). 
Let $W \dasharrow Y_1 \to X$ be the end result. 
Since 
\[
K_W + \widetilde{\Delta} ^{\wedge 1} + F \sim_{g, \R}  
-((\widetilde{\Delta}-\widetilde{\Delta} ^{\wedge 1})+E-F),
\]
this MMP is also a 
$- ((\widetilde{\Delta} - \widetilde{\Delta} ^{\wedge 1}) + E - F)$-MMP over $X$. 
Hence by the negativity lemma, the push-forward of $(\widetilde{\Delta} - \widetilde{\Delta} ^{\wedge 1}) + E - F$ 
on the end result $Y_1$ is effective. 
This implies that the $g$-exceptional divisors $E_i$ with $\mathrm{coeff}_{E_i} E < 1$ are contracted in this MMP. 
Therefore the condition (1) holds. 
Since $\Delta _{Y_1}  ^{\wedge 1}$ is nothing but the push-forward of $\widetilde{\Delta} ^{\wedge 1} + F$, 
the condition (2) follows because being dlt is preserved under an MMP. 
This completes the proof of Step \ref{s1-dlt-modif}. 
\end{proof}

\begin{step}\label{s2-dlt-modif} 
There exists a projective birational morphism $f:Y \to X$ that satisfies 
the conditions 
(1)-(3) in the statement. 
\end{step}

\begin{proof}[Proof of Step \ref{s2-dlt-modif}] 
The following argument is very similar to \cite[Lemma 29]{dFKX}. 
By Step \ref{s1-dlt-modif}, 
we can find a projective birational morphism $f_1 : Y_1 \to X$  with the conditions (1) and (2). 
Define $\Delta _{Y_1}$ by $f_1 ^* (K_X + \Delta) = K_{Y_1} + \Delta _{Y_1}$. 

There exists a $(K_{Y_1} + \Delta _{Y_1} ^{<1})$-MMP over $X$ that terminates (Theorem \ref{t-lc-eff-mmp}). 
Let $Y_1 \dasharrow Y_2 \overset{f_2}{\longrightarrow} X$ be the end result, and 
let $\Delta _{Y_2}$ be the push-forward of $\Delta_{Y_1}$ on $Y_2$. 
We obtain 
\begin{equation}\label{e1-dlt-modif} 
\Supp (\Delta _{Y_2} ^{\ge 1}) = \Nklt (Y_2, \Delta _{Y_2}), 
\end{equation}
since $(Y_2, \Delta _{Y_2} ^{<1})$ is klt. 

Let us show the equation 
\begin{equation}\label{e2-dlt-modif} 
\Nklt (Y_2, \Delta _{Y_2}) = f_2 ^{-1} (\Nklt (X, \Delta)). 
\end{equation}
By $f_2 ^{*} (K_X + \Delta) = K_{Y_2} + \Delta_{Y_2}$, 
$f_2$ induces a surjective morphism $\Nklt (Y_2, \Delta _{Y_2}) \to \Nklt (X, \Delta)$ 
(cf.\ \cite[Lemma 2.30]{KM98}). 
In particular, the inclusion $\Nklt (Y_2, \Delta _{Y_2}) \subset f_2 ^{-1} (\Nklt (X, \Delta))$ 
holds. 
Assuming that $\Nklt (Y_2, \Delta _{Y_2}) \not\supset f_2 ^{-1} (\Nklt (X, \Delta))$, 
let us derive a contradiction. 
Take a closed point $y \in f_2 ^{-1} (\Nklt (X, \Delta)) \setminus \Nklt (Y_2, \Delta _{Y_2})$. 
Set $x:=f_2(y)$. 
Then it holds that $x \in \Nklt (X, \Delta)$ and $f_2 ^{-1}(x) \not\subset \Nklt (Y_2, \Delta _{Y_2})$. 
Since $\Nklt (Y_2, \Delta _{Y_2}) \to \Nklt (X, \Delta)$ is surjective, 
there exists a curve $C$, contained in $f^{-1}_{2}(x)$, such that 
$C$ intersects but is not contained in $\Nklt (Y_2, \Delta _{Y_2})$. 
By (\ref{e1-dlt-modif}), it holds that $\Delta _{Y_2} ^{\ge 1} \cdot C>0$. 
On the other hand, 
\[
-\Delta _{Y_2} ^{\geq 1} = - (\Delta _{Y_2} - \Delta _{Y_2} ^{<1} ) 
\equiv_{f_2} K_{Y_2}+\Delta_{Y_2}^{< 1} 
\] is $f_2$-nef. 
This is a contradiction. 
Hence, (\ref{e2-dlt-modif}) holds.

Let $h: (Y, \Delta _Y) \to (Y_2, \Delta _{Y_2})$ be a projective birational morphism 
satisfying (1) and (2), whose existence is guaranteed by Step 1. 
Since $Y_2$ is $\mathbb{Q}$-factorial, $\mathrm{Ex}(h)$ is purely codimension one, and hence
\[
\Nklt (Y, \Delta _{Y}) = h ^{-1} (\Nklt (Y_2, \Delta _{Y_2}))
\]
holds. 
Therefore the composition $f:=f_2 \circ h : Y \to X$ satisfies the conditions (1)-(3). 
This completes the proof of Step \ref{s2-dlt-modif}. 
\end{proof}

Step \ref{s2-dlt-modif} completes the proof of Proposition \ref{p-dlt-modif}. 
\end{proof}

\begin{cor}\label{c-dlt-modif} 
Let $k$ be a perfect field of characteristic $p>5$. 
Let $(X, \Delta)$ be a three-dimensional quasi-projective log canonical pair over $k$. 
Then there exists a projective birational $k$-morphism $f:Y \to X$ 
from a normal threefold $Y$ over $k$ that satisfies the following conditions: 
\begin{enumerate}
\item $a_F (X, \Delta) = 0$ holds for any $f$-exceptional prime divisor $F$.  
\item $(Y, \Delta _Y)$ is a $\Q$-factorial dlt pair,  
where $\Delta _Y$ is the $\R$-divisor defined  by $K_Y + \Delta _Y = f^* (K_X + \Delta)$. 
\item 
$\Nklt (Y, \Delta _{Y}) = f ^{-1} (\Nklt (X, \Delta))$ holds. 
\end{enumerate} 
\end{cor}

\begin{proof}
The claim directly follows from Proposition \ref{p-dlt-modif}. 
\end{proof}

\section{Cone theorem}

The purpose of this section is to prove the cone theorem for log canonical threefolds 
(Theorem \ref{t-lc-cone}, Theorem \ref{t-lc-cone2}). 
To this end, we start with a projective case (Subsection \ref{ss1-cone-thm}). 
Then we treat the dlt case (Subsection \ref{ss2-cone-thm}). 
Finally, we will establish the cone theorem for log canonical threefolds (Subsection \ref{ss3-cone-thm}). 
As a consequence of the cone theorem, 
we obtain a result on the Shokurov polytope (Subsection \ref{ss4-cone-thm}).

\subsection{Projective case}\label{ss1-cone-thm}

In this subsection, we prove Lemma \ref{l-proj-cone}. 
Let us start with a criterion to deduce the cone theorem. 

\begin{lem}\label{l-cone-criterion}
Let $k$ be a field. 
Let $f:X \to Z$ be a projective $k$-morphism 
from a normal $k$-variety $X$ to a quasi-projective $k$-scheme $Z$. 
Let $\Delta$ be an $\R$-divisor on $X$ such that $K_X+\Delta$ is $\R$-Cartier. 
Let $A$ be an $f$-ample $\R$-Cartier $\R$-divisor on $X$. 
For any ample $\R$-Cartier $\R$-divisor $H$, 
set $a_H$ to be the $f$-nef threshold of $(K_X+\Delta+\frac{1}{2}A, H)$. 
Assume that there exist finitely many curves $C_1, \ldots, C_m$ on $X$ 
such that 
\begin{enumerate}
\item 
$f(C_i)$ is a point for any $i\in \{1, \ldots, m\}$, and 
\item 
for any ample $\R$-Cartier $\R$-divisor $H$ on $X$, it holds that 
$(K_X+\Delta+\frac{1}{2}A+a_H H)\cdot C_i=0$ for some $i \in \{1, \ldots, m\}$. 
\end{enumerate}
Then the following equation holds: 
\[
\overline{\NE}(X/Z)=\overline{\NE}(X/Z)_{K_X+\Delta+A \geq 0} +\sum_{i=1}^m \R_{\geq 0}[C_i]. 
\]
\end{lem}

\begin{proof}
We can apply the same argument as in \cite[Lemma 6.2]{CTX15}. 
\end{proof}

\begin{lem}\label{l-proj-cone}
Let $k$ be an algebraically closed field of characteristic $p>5$. 
Let $(X, \Delta)$ be a three-dimensional log canonical pair over $k$ and 
let $f:X \to Z$ be a projective $k$-morphism to a projective $k$-scheme $Z$. 
Let $A$ be an $f$-ample $\R$-Cartier $\R$-divisor on $X$. 
Then there exist finitely many curves $C_1, \ldots, C_m$ on $X$ such that 
\begin{enumerate}
\item 
$f(C_i)$ is a point for any $i\in \{1, \ldots, m\}$, and  
\item 
$\overline{\NE}(X/Z)=\overline{\NE}(X/Z)_{K_X+\Delta+A \geq 0} +\sum_{i=1}^m \R_{\geq 0}[C_i]. $
\end{enumerate}
\end{lem}

\begin{proof}
Taking the Stein factorisation of $f$, we may assume that $f_*\MO_X=\MO_Z$. 
Replacing $A$ by $A+f^*A_Z$ for some ample Cartier divisor $A_Z$ on $Z$, 
the problem is reduced to the case when $A$ is ample. 
Let $H_Z$ be an ample Cartier divisor on $Z$. 
By \cite[Theorem 1]{Tanb}, 
there exists an effective $\R$-divisor $M$ such that $M \sim_{\R} 7f^*H_Z$ 
and $(X, \Delta':=\Delta+M)$ is log canonical. 
Applying Theorem \ref{t-Wal1.7} to $(X, \Delta)$ and $\frac{1}{2}A$, 
there exist finitely many curves $C_1, \ldots, C_n$ with the condition 
$0 > \left( K_X + \Delta + \frac{1}{2} A \right) \cdot C_i \ge -6$ and  
\begin{equation}\label{e-proj-cone}
\overline{\NE}(X)=\overline{\NE}(X)_{K_X+\Delta+\frac{1}{2}A \geq 0} +\sum_{i=1}^n \R_{\geq 0}[C_i]. 
\end{equation}
After permuting indices, we can find an integer $m$ such that 
\begin{itemize}
\item $0 \leq m \leq n$,
\item $f(C_{\alpha})$ is a point for any $\alpha \in \{1, \ldots, m\}$, and 
\item $f(C_{\beta})$ is not a point for any $\beta \in \{m+1, \ldots, n\}$. 
\end{itemize}
It is enough to prove that the curves $C_1, \ldots, C_m$ satisfy 
the condition (2). 
Take an ample $\R$-Cartier $\R$-divisor $H$ on $X$ and 
let $a_H$ be the $f$-nef threshold of $(K_X+\Delta'+\frac{1}{2}A, H)$. 
Note that $a_H$ is equal to the nef threshold of $(K_X+\Delta'+\frac{1}{2}A, H)$. 
Indeed, for any $\beta \in \{m+1, \ldots, n\}$, it holds that 
\[
\big( K_X+\Delta'+\frac{1}{2}A + a_H H \big) \cdot C_{\beta}
\ge \big( K_X+\Delta+7f^*H_Z +\frac{1}{2}A \big)\cdot C_{\beta} 
\geq -6+7=1. 
\]
Therefore, by (\ref{e-proj-cone}), 
there exists $i \in \{1, \ldots, n\}$ such that 
$(K_X+\Delta'+\frac{1}{2}A+a_H H)\cdot C_i=0$. 
By the inequality above, it holds that $1 \leq i \leq m$, as desired. 
\end{proof}

\subsection{Dlt case}\label{ss2-cone-thm}

The purpose of this subsection is to prove the cone theorem for the dlt case (Proposition \ref{p-dlt-cone}). 
To this end, we first find an extremal ray that can be compactified 
(Lemma \ref{l-dlt-rationality}). 
We also need a basic fact on extremal rays (Lemma \ref{l-ext-length}).

\begin{lem}\label{l-dlt-rationality}
Let $k$ be an algebraically closed field of characteristic $p>5$. 
Let $(X, \Delta)$ be a three-dimensional $\Q$-factorial log canonical pair over $k$ such that $(X, \{\Delta\})$ is klt. 
Let $f:X \to Z$ be a projective $k$-morphism 
to a quasi-projective $k$-scheme $Z$. 
Let $A$ be an $f$-ample $\R$-divisor on $X$ 
such that $K_X+\Delta+A$ is $f$-nef but not $f$-ample. 
Then there exists a commutative diagram of quasi-projective $k$-schemes: 
$$\begin{CD}
X @>j_X>> \overline X\\
@VVfV @VV\overline fV\\
Z @>j_Z>> \overline Z.
\end{CD}$$
such that 
\begin{enumerate}
\item[(a)] 
$j_X$ and $j_Z$ are open immersions to projective $k$-schemes 
$\overline{X}$ and $\overline{Z}$, 
\item[(b)] 
$(\overline X, \overline{\Delta})$ is a $\Q$-factorial log canonical pair 
such that $(\overline X, \{\overline \Delta\})$ is klt, where $\overline \Delta$ denotes the closure of $\Delta$ in $\overline X$, and 
\item[(c)]  
there is a $(K_{\overline X}+\overline \Delta)$-negative extremal ray $R$ of 
$\overline{\NE}(\overline X/\overline Z)$ such that 
for the contraction $\overline \psi:\overline X \to \overline Y$ of $R$, 
its restriction $\psi:X \to Y$ over $Z$ is not an isomorphism and 
$K_X+\Delta+A$ is $\psi$-numerically trivial. 
\end{enumerate}
\end{lem}

\begin{proof}
Fix an open immersion $j_Z:Z \to \overline Z$ to a projective $k$-scheme $\overline Z$. 
We divide the proof into three steps. 

\setcounter{step}{0}

\begin{step}\label{s1-dlt-rationality}
Assume that $(X, \Delta)$ is klt. 
Then there exist an open immersion $j^{(1)}:X \hookrightarrow X^{(1)}$ to a $\mathbb{Q}$-factorial projective threefold $X^{(1)}$, a projective morphism 
$f^{(1)}:X^{(1)} \to \overline Z$ and effective $\R$-divisors $\Delta^{(1)}$ and $A^{(1)}$ on $X^{(1)}$ such that 
\begin{enumerate}
\item
$j_Z \circ f=f^{(1)}\circ j^{(1)}$,  
\item
$\Delta^{(1)}$ is the closure of $\Delta$, 
\item 
$A^{(1)}|_X \sim_{Z, \R} A$, 
\item 
$(X^{(1)}, \Delta^{(1)}+A^{(1)})$ is klt, and  
\item 
$K_{X^{(1)}}+ \Delta^{(1)}+A^{(1)}$ is $f^{(1)}$-nef. 
\end{enumerate}
\end{step}

\begin{proof}[Proof of Step \ref{s1-dlt-rationality}]
By Lemma \ref{l-lcklt-perturb} (2), 
we may assume that $A$ is effective and $(X, \Delta+A)$ is klt. 
By Proposition \ref{p-klt-cpt}, there is an open immersion $j:X \to  X^{(0)}$ over $\overline Z$ to a normal $\Q$-factorial threefold $X^{(0)}$ projective over $\overline Z$ such that 
$(X^{(0)}, \Delta^{(0)}+A^{(0)})$ is klt, where $\Delta^{(0)}$ and $A^{(0)}$ denote the closures of $\Delta$ and $A$ respectively. 

It follows from \cite[Theorem 1.6]{BW17} that 
there exists a $(K_{X^{(0)}}+\Delta^{(0)}+A^{(0)})$-MMP over $\overline Z$ 
that terminates with a log minimal model $(X^{(1)}, \Delta^{(1)}+A^{(1)})$. 
Since $(K_{X^{(0)}}+\Delta^{(0)}+A^{(0)})|_X=K_X+\Delta+A$ is nef over $Z$, 
the restriction of $X^{(0)} \dashrightarrow X^{(1)}$ over $Z$ is an isomorphism. 
By construction, all the properties (1)-(5) hold. 
This completes the proof of Step \ref{s1-dlt-rationality}. 
\end{proof}

\begin{step}\label{s2-dlt-rationality}
The assertion of Lemma \ref{l-dlt-rationality} holds if $(X, \Delta)$ is klt and $\Delta$ is big over $Z$. 
\end{step}

\begin{proof}[Proof of Step \ref{s2-dlt-rationality}]
Take $j^{(1)}:X \hookrightarrow X^{(1)}$, $f^{(1)}:X^{(1)} \to \overline Z$, $\Delta^{(1)}$ and $A^{(1)}$ as in Step \ref{s1-dlt-rationality}. 
By \cite[Theorem 1.5]{BW17}, 
there exists a $(K_{X^{(1)}}+\Delta^{(1)})$-MMP over $\overline Z$ with scaling of $A^{(1)}$ 
that terminates: 
\[
X^{(1)}=:X^{(1)}_0 \overset{\varphi_0}{\dashrightarrow} X^{(1)}_1 \overset{\varphi_1}{\dashrightarrow} \cdots \overset{\varphi_{\ell-1}}{\dashrightarrow} X^{(1)}_{\ell}. 
\]
For any $i$, we denote the associated morphism 
corresponding to the extremal ray by  $\overline{\psi}_{i}:X^{(1)}_i\to \overline{Y}_{i}$. 
Then there are two possibilities as follows.  
\begin{enumerate}
\item[(A)] For any $i \in \{0, \ldots, \ell - 1 \}$, the restriction of $\varphi_i$ over $Z$ is an isomorphism. 
\item[(B)] There exists $i \in \{0, \ldots, \ell - 1 \}$ such that the restriction of $\varphi_i$ over $Z$ is not an isomorphism. 
\end{enumerate}

Assume that (A) occurs. 
Then we set $(\overline X, \overline \Delta):=(X^{(1)}_{\ell}, \Delta^{(1)}_{\ell})$. 
It follows from the construction that (a) and (b) hold. 
In particular, $K_{\overline X}+\overline \Delta$ is not nef over $\overline Z$, 
since its restriction $(K_{\overline X}+\overline \Delta)|_X=K_X+\Delta$ is not nef over $Z$. 
Hence, there is a $(K_{\overline X}+\overline \Delta)$-negative extremal ray $R$ of 
$\overline{\NE}(\overline X/\overline Z)$ that corresponds to 
a $(K_{\overline X}+\Delta)$-Mori fibre space $\overline \psi:\overline X \to \overline Y$ 
over $\overline Z$. 
It is clear that the restriction of $\overline \psi$ over $Z$ is not an isomorphism, 
hence (c) holds.

Assume that (B) occurs. 
Then there is an index $j \in \{0, \ldots, \ell\}$ such that the restriction of $\varphi_j$ over $Z$ is not an isomorphism and 
the restriction of $\varphi_{j'}$ over $Z$ is an isomorphism over $Z$ for any $j'<j$. 
In this case, we set $(\overline X, \overline \Delta):=(X^{(1)}_j, \Delta^{(1)}_j)$, 
$\overline Y:=\overline Y_j$ and $\overline \psi:=\overline{\psi}_{j}$. 
Then all the properties (a), (b) and (c) hold by construction. 
This completes the proof of Step \ref{s2-dlt-rationality}. 
\end{proof}

\begin{step}\label{s3-dlt-rationality}
The assertion of Lemma \ref{l-dlt-rationality} holds without any additional assumptions. 
\end{step}

\begin{proof}[Proof of Step \ref{s3-dlt-rationality}]
Pick a sufficiently small positive real number $\epsilon$ such that $\frac{1}{2}A+\epsilon \lfloor \Delta\rfloor$ is $f$-ample. 
Since $(X,\Delta-\epsilon \lfloor \Delta\rfloor)$ is klt, there exists an effective $\R$-divisor 
$A' \sim_{f, \R}\frac{1}{2}A+\epsilon \lfloor \Delta\rfloor$ 
such that the pair $(X, \Delta':=\Delta-\epsilon \lfloor \Delta\rfloor+A')$ is klt by Lemma \ref{l-lcklt-perturb} (2). 
Again by Lemma \ref{l-lcklt-perturb} (2), there is an effective $\R$-divisor $A''$ on $X$ such that $A'' \sim_{f, \R}\frac{1}{2}A$ and $(X,\Delta'+A'')$ is klt. 
By construction, $\Delta'$ is $f$-big and $K_{X}+\Delta+A\sim_{f, \R}K_{X}+\Delta'+A''$. 
In particular, $K_{X}+\Delta'+A''$ is $f$-nef but not $f$-ample. 

Applying Step \ref{s2-dlt-rationality} to $f:X \to Z$, $(X, \Delta')$ and $A''$, 
there exists a commutative diagram of quasi-projective $k$-schemes: 
$$\begin{CD}
X @>j_X>> \overline X\\
@VVfV @VV\overline fV\\
Z @>j_Z>> \overline Z.
\end{CD}$$
such that 
\begin{enumerate}
\item[(a)] 
$j_X$ and $j_Z$ are open immersions to projective $k$-schemes 
$\overline{X}$ and $\overline{Z}$, 
\item[(b)'] 
$(\overline X, \overline{\Delta'})$ is a $\Q$-factorial klt pair, 
where $\overline{\Delta'}$ denotes the closure of $\Delta'$ in $\overline X$, and 
\item[(c)']  
there is a $(K_{\overline X}+\overline{\Delta'})$-negative extremal ray $R$ of 
$\overline{\NE}(\overline X/\overline Z)$ such that 
for the contraction $\overline \psi:\overline X \to \overline Y$ of $R$, 
its restriction $\psi:X \to Y$ over $Z$ is not an isomorphism and 
$K_X+\Delta'+A''$ is $\psi$-numerically trivial. 
\end{enumerate}

It suffices to prove that (b) and (c) hold. 
We first show (b). 
By $\Delta'=\Delta-\epsilon \lfloor \Delta\rfloor +A'$, 
we obtain $\overline{\Delta'} \geq \overline{\Delta}-\epsilon \lfloor \overline{\Delta}\rfloor$.
Since $\overline X$ is $\Q$-factorial and the log pair $(\overline X, \overline{\Delta'})$ is klt, we see that $(\overline X, \overline{\Delta}-\epsilon \lfloor \overline{\Delta}\rfloor)$ is klt. 
Hence, also $(\overline X, \{\overline{\Delta}\})$ is klt. 
Furthermore, the log canonical threshold ${\rm lct}(\overline X, \{\overline{\Delta}\};\lfloor \overline{\Delta}\rfloor)$ is greater than $1-\epsilon$. 
As $\epsilon>0$ was chosen to be sufficiently small, 
ACC for log canonical thresholds implies that $(\overline X, \overline \Delta)$ is log canonical. 
Thus, (b) holds. 

Let us show (c). 
Pick a curve $C$ on $X$ contracted by $\psi$, 
whose existence is guaranteed by (c)'. 
Then $C$ generates $R$ and therefore we have $(K_{X}+ \Delta'+A'')\cdot C=0$ and $(K_{X}+ \Delta')\cdot C<0$. 
Now recall that $K_{X}+\Delta+A\sim_{f,\R}K_{X}+\Delta'+A''$ and $K_{X}+\Delta+\frac{1}{2}A\sim_{f,\R}K_{X}+\Delta'$, which imply $(K_{X}+ \Delta+A)\cdot C=0$ and $(K_{X}+ \Delta)\cdot C<0$. 
Thus we see that $K_{X}+ \Delta+A$ is $\psi$-numerically trivial and 
$R$ is a $(K_{\overline X}+\overline{\Delta})$-negative extremal ray of $\overline{\NE}(\overline X/\overline Z)$. 
Hence, (c) holds. 
This completes the proof of Step \ref{s3-dlt-rationality}. 
\end{proof}
Step \ref{s3-dlt-rationality} completes the proof of Lemma \ref{l-dlt-rationality}. 
\end{proof}

\begin{lem}\label{l-ext-length}
Let $k$ be an algebraically closed field of characteristic $p>5$. 
Let $(X, \Delta)$ be a three-dimensional $\Q$-factorial log canonical pair over $k$ 
and let $f:X \to Z$ be a projective $k$-morphism to a projective $k$-scheme. 
Let $R$ be a $(K_X+\Delta)$-negative extremal ray of $\overline{\NE}(X/Z)$ and 
let $\varphi:X \to Y$ be the contraction of $R$. 
\begin{enumerate}
\item 
If $\varphi$ is birational, 
then there exists a non-empty open subset $Y'$ of $\varphi(\Ex(\varphi))$ 
such that for any closed point $y \in Y'$, 
there is a rational curve $C$ on $X$ such that $\varphi(C)=\{y\}$ and 
$0<-(K_X+\Delta) \cdot C \leq 6$. 
\item 
If $\varphi$ is not birational, 
then 
for any closed point $y \in Y$, 
there is a rational curve $C$ on $X$ such that $\varphi(C)=\{y\}$ and 
$0<-(K_X+\Delta) \cdot C \leq 6$. 
\end{enumerate}
\end{lem}

\begin{proof}
We may assume that $Z=\Spec\,k$. 
By the cone theorem, 
we can find an ample $\R$-divisor $H$ such that 
$\overline{\NE}(X) \cap (K_X+\Delta+H)^{\perp}=R$. 
If $K_X+\Delta+H$ is not big i.e.\ $\varphi$ is not birational, 
then the assertion follows from \cite[Corollary 1.5]{CTX15}. 
Hence, we may assume that $\varphi$ is birational. 
If $\dim\Ex(\varphi)=1$, then the assertion holds by 
Theorem \ref{t-Wal1.7}. 
Therefore, the problem is reduced to the case when 
$\varphi$ is a divisorial contraction. 
In this case, we can apply the same argument as in 
\cite[the second paragraph of the proof of Lemma 3.2]{BW17} 
without any changes. 
\end{proof}

\begin{prop}\label{p-dlt-cone}
Let $k$ be an algebraically closed field of characteristic $p>5$. 
Let $(X, \Delta)$ be a three-dimensional $\Q$-factorial log canonical pair over $k$ such that $(X, \{\Delta\})$ is klt. 
Let $f:X \to Z$ be a projective $k$-morphism to a quasi-projective $k$-scheme $Z$. 
Then the following hold. 
\begin{enumerate}
\item 
If $H$ is an ample $\R$-divisor on $X$ such that $K_X+\Delta+H$ is $f$-nef but not $f$-ample, 
then there exists a rational curve $C$ such that $f(C)$ is a point, $(K_X+\Delta+H) \cdot C = 0$ and 
$0<-(K_X+\Delta) \cdot C \leq 6$. 

\item 
If $A$ is an $f$-ample $\R$-divisor on $X$, then 
there exist finitely many rational curves $C_1, \ldots, C_m$ such that 
\[
\overline{\NE}(X/Z)=\overline{\NE}(X/Z)_{K_X+\Delta+A \geq 0}+\sum_{i=1}^m \R_{\geq 0} [C_i].
\]
\end{enumerate}
\end{prop}

\begin{proof}
The assertion (1) follows from Lemma \ref{l-dlt-rationality} and Lemma \ref{l-ext-length}. 
We prove (2). 
Let $H$ be an $f$-ample $\R$-divisor on $X$ and let $a_H$ be the $f$-nef threshold of $(K_X+\Delta+\frac{1}{2}A, H)$. 
By construction $K_X+\Delta+(\frac{1}{2}A+a_H H)$ is $f$-nef but not $f$-ample. 
By (1), there exists a rational curve $C_H$ such that $f(C_H)$ is a point, $(K_X+\Delta+\frac{1}{2}A+a_H H)\cdot C_H=0$ and $0<-(K_X+\Delta) \cdot C_H \leq 6$. 
From this we have 
\[
0<A \cdot C_H \leq (A+2a_H H)\cdot C_H=-2(K_X+\Delta)\cdot C_H\leq6\cdot 2=12.
\]
Therefore the subset $\{[C_H]\in \overline{\NE}(X/Z)| \,H\,{\rm is}\,f\mathchar`-{\rm ample}\}$ of $\overline{\NE}(X/Z)$ is a finite set. 
In this way we can find finitely many rational curves $C_{H_{1}},\ldots, C_{H_{m}}$ satisfying the conditions of Lemma \ref{l-cone-criterion}. Hence, (2) follows.  
\end{proof}

\subsection{Log canonical case}\label{ss3-cone-thm}

Let us prove the main results in this section (Theorem \ref{t-lc-cone}, Theorem \ref{t-lc-cone2}). 

\begin{thm}\label{t-lc-cone}
Let $k$ be an algebraically closed field of characteristic $p>5$. 
Let $(X, \Delta)$ be a three-dimensional log canonical pair over $k$ 
and let $f:X \to Z$ be a projective $k$-morphism to a quasi-projective $k$-scheme $Z$. 
Then there exists a countable set $\{ C_i \} _{i \in I}$ 
of rational curves on $X$ which satisfies the following conditions: 
\begin{enumerate}
\item $f(C_i)$ is a point for any $i \in I$. 
\item $\overline{\NE}(X/Z)=\overline{\NE}(X/Z)_{K_X+\Delta \geq 0}+\sum_{i \in I} \R_{\geq 0}  [C_i]$. 
\item $0 < - (K_X + \Delta) \cdot C_i \le 6$ for any $i \in I$. 
\item 
For any $f$-ample $\R$-Cartier $\mathbb{R}$-divisor $A$, 
there exists a finite subset $J$ of $I$ such that 
\[
\overline{\NE}(X/Z)=\overline{\NE}(X/Z)_{K_X+\Delta+A \geq 0}+\sum_{j \in J} \R_{\geq 0}  [C_j].
\] 
\end{enumerate}
\end{thm}

\begin{proof}
The proof consists of two steps. 
\setcounter{step}{0}
\begin{step}\label{s1-lc-cone}
If $X$ is $\Q$-factorial and $(X, \{\Delta\})$ is klt, 
then there exists a countable set $\{ C_i \} _{i \in I}$ 
of rational curves on $X$ that satisfies the conditions (1)--(3) of Theorem \ref{t-lc-cone}. 
\end{step}
\begin{proof}
Take any $f$-ample $\R$-Cartier $\mathbb{R}$-divisor $H$. Then for any positive integer $n$, 
Proposition \ref{p-dlt-cone} enables us to find 
a finite set of rational curves $\{ C_i \} _{i \in I_n}$ such that 
$f(C_i)$ is a point, $-(K_X+\Delta) \cdot C_i \leq 6$, and 
\[
\overline{\NE}(X/Z)=\overline{\NE}(X/Z)_{K_X+\Delta + \frac{1}{n}H \geq 0}+\sum_{i \in I_n} \R_{\geq 0}  [C_i]. 
\]
Therefore $I = \bigcup _{n \ge 1} I_n$ satisfies the conditions (1)--(3). 
This completes the proof of Step \ref{s1-lc-cone}. 
\end{proof}

\begin{step}\label{s2-lc-cone}
The assertion of Theorem \ref{t-lc-cone} holds without any additional assumptions. 
\end{step}

\begin{proof}
Let $g: Y \to X$ be a projective birational morphism with the conditions in Corollary \ref{c-dlt-modif}. 
We apply Step \ref{s1-lc-cone} to $(Y, \Delta _Y)$ to conclude that there exists 
a countable set $\{ C_{i}^{Y} \} _{i \in J}$
of rational curves on $Y$ with the following conditions: 
\begin{itemize}
\item $f(g(C_{i}^{Y}))$ is a point for any $i \in J$. 
\item $\overline{\NE}(Y/Z)=\overline{\NE}(Y/Z)_{K_Y + \Delta _Y \geq 0}+\sum_{i \in J} \R_{\geq 0} [C_i ^Y]$. 

\item $0 < - (K_Y + \Delta _Y) \cdot C_i ^Y \le 6$ for each $i \in J$. 
\end{itemize}
Set $C_i:=g(C_i ^Y)$ and $I := \{ i \in J \mid \dim C_i = 1 \}$. 
Then it is clear that (1) holds. 
We get (2) by $\overline{\NE}(X/Z) = g_* (\overline{\NE}(Y/Z))$. 
It follows from the projection formula that (3) holds. 

We now show (4). 
Let $R$ be an extremal ray of $\overline{\NE}(X/Z)$ which is $(K_X + \Delta + A)$-negative. 
By (1) and (2) we have already proved, we obtain $R=\R_{\geq 0}[C]$ 
for some curve $C$ such that $f(C)$ is a point and $0 < - (K_X + \Delta) \cdot C \le 6$. 
Then it holds that 
\[
A \cdot C <-(K_X+\Delta) \cdot C \leq 6. 
\]
Hence, there are only finitely many extremal rays satisfying this property (cf.\ \cite[Corollary 1.19 (3)]{KM98}). 
Thus (4) holds. 
This completes the proof of Step \ref{s2-lc-cone}. 
\end{proof}
Step \ref{s2-lc-cone} completes the proof of  Theorem \ref{t-lc-cone}. 
\end{proof}

\begin{thm}\label{t-lc-cone2}
Let $k$ be a perfect field of characteristic $p>5$. 
Let $(X, \Delta)$ be a three-dimensional log canonical pair over $k$ 
and let $f:X \to Z$ be a projective $k$-morphism to a quasi-projective $k$-scheme $Z$. 
Let $A$ be an $f$-ample $\R$-Cartier $\R$-divisor on $X$. 
Then there exist finitely many curves $C_1, \ldots, C_m$ on $X$ such that 
\begin{enumerate}
\item $f(C_i)$ is a point for any $i \in \{1, \ldots, m\}$, and 
\item 
$\overline{\NE}(X/Z)=\overline{\NE}(X/Z)_{K_X+\Delta+A \geq 0}+\sum_{i=1}^m \R_{\geq 0}  [C_i]$.
\end{enumerate}
\end{thm}

\begin{proof}
We may assume that 
\begin{itemize}
\item $f_*\MO_X=\MO_Z$, 
\item $k$ is algebraically closed in $K(Z)$, and 
\item both $X$ and $Z$ are geometrically integral and geometrically normal over $k$. 
\end{itemize}
Indeed, we may assume the first condition by taking the Stein factorisation of $f$. 
Then, after replacing $k$ by the algebraic closure of $k$ in $K(Z)$, 
the second condition holds. 
The third condition automatically follows from the other two. 

Let $\overline k$ be the algebraic closure of $k$. 
We set  
\[
f_{\overline k}:X_{\overline k}\to Z_{\overline k}
\]
to be the base change $f \times_k \overline k$. 
Let $\Delta_{\overline k}$ and $A_{\overline k}$ be the $\R$-divisors defined as 
the pullbacks of $\Delta$ and $A$ respectively. 
Then $(X_{\overline k}, \Delta_{\overline k}, f_{\overline k}, A_{\overline k})$ 
satisfies the assumptions listed in the statement. 
By Theorem \ref{t-lc-cone}, 
there exist finitely many curves $C'_1, \ldots, C'_m$ on $X_{\overline k}$ such that 
\begin{enumerate}
\item[(i)] $f_{\overline k}(C'_i)$ is a point for any $i$, and 
\item[(ii)] 
$\overline{\NE}(X_{\overline k}/Z_{\overline k})=\overline{\NE}(X_{\overline k}/Z_{\overline k})_{K_{X_{\overline k}}+\Delta_{\overline k}+\frac{1}{2}A_{\overline k} \geq 0}+\sum_{i=1}^m \R_{\geq 0} [C'_i]$. 
\end{enumerate}
Let $C_i$ be the image of $C'_i$ for any $i \in \{1, \ldots, m\}$. 
It follows from (i) that the condition Lemma \ref{l-cone-criterion}(1) holds. 
Therefore, it is enough to prove that also Lemma \ref{l-cone-criterion}(2) holds. 
Let $H$ be an ample $\R$-Cartier $\R$-divisor on $X$ and 
let $a_H$ be the $f$-nef threshold of $(K_X+\Delta+\frac{1}{2} A, H)$. 
We see that $a_H$ is equal to the $f_{\overline{k}}$-nef threshold of $(K_{X_{\overline k}}+\Delta_{\overline k}+\frac{1}{2} A_{\overline k}, H_{\overline k})$ 
(cf.\ \cite[Lemma 2.3]{Tana}). 
Hence, by (ii), we get 
$(K_{X_{\overline k}}+\Delta_{\overline k}+\frac{1}{2} A_{\overline k}+a_H H_{\overline k}) \cdot C'_i=0$ for some $i \in \{1,\ldots, m\}$, 
which in turn implies  $(K_X+\Delta+\frac{1}{2} A+a_H H) \cdot C_i=0$ (cf.\ \cite[Lemma 2.3]{Tana}). 
Therefore, Lemma \ref{l-cone-criterion}(2) holds, as desired. 
\end{proof}

\subsection{Shokurov polytope}\label{ss4-cone-thm}

As a consequence of the cone theorem, we obtain a result on the Shokurov polytope. 
We first fix some terminologies. 

\begin{nota}\label{n-polytope}
Let $k$ be an algebraically closed field of characteristic $p>5$. 
Let $X$ be a $\Q$-factorial klt threefold and let $f:X \to Z$ be a projective morphism 
to a quasi-projective $k$-scheme. 
Fix finitely many prime divisors $D_1, \ldots, D_n$ and set 
\[
V:=\bigoplus_{i=1}^n \R \cdot D_i,
\]
which is a subspace of the $\R$-vector space of the $\R$-divisors on $X$. 
For any $D=\sum d_i D_i \in V$, we set $||D||:=\max_{1 \leq i \leq n} \{|d_i|\}$. 
We see that 
\[
\mathcal L:=\{\Delta \in V\,|\, (X, \Delta)\text{ is log canonical}\}
\]
is a rational polytope in $V$. 
\end{nota}

\begin{prop}\label{p-polytope}
We use Notation \ref{n-polytope}. 
Fix $D \in \mathcal L$. 
Then there exist positive real numbers $\alpha$ and $\delta$ which satisfy the following properties. 
\begin{enumerate}
\item 
If $\Gamma$ is an extremal curve of $\overline{\NE}(X/Z)$ 
and if $(K_X+D) \cdot \Gamma>0$, 
then $(K_X+D) \cdot \Gamma>\alpha$. 
\item 
If $\Delta \in \mathcal L$, $||\Delta-D||<\delta$, and $(K_X+\Delta) \cdot R \leq 0$ 
for an extremal ray $R$ of $\overline{\NE}(X/Z)$, 
then $(K_X+D) \cdot R \leq 0$. 
\item 
Let $\{R_t\}_{t \in T}$ be a set of extremal rays of $\overline{\NE}(X/Z)$. 
Then the set 
\[
\mathcal N_T:=\{\Delta \in \mathcal L\,|\,(K_X+\Delta) \cdot R_t\geq 0 
\text{ for any }t\in T\}
\]
is a rational polytope. 
\item 
Assume that $K_X+D$ is $f$-nef, $\Delta \in \mathcal L$ and that 
\[
X=X_1 \dashrightarrow X_2 \dashrightarrow \cdots 
\]
is a sequence of $(K_X+\Delta)$-MMP over $Z$ that consists of flips which are 
$(K_X+D)$-numerically trivial. 
Then, for any $i$ and any curve $\Gamma$ on $X_i$ whose image on $Z$ is a point, 
if $(K_{X_i}+D_i) \cdot \Gamma >0$, then $(K_{X_i}+D_i) \cdot \Gamma>\alpha$, 
where $D_i$ denotes the push-forward of $D$ on $X_i$. 
\item 
In addition to the assumptions of (4), suppose that $||\Delta-D||<\delta$. 
If $(K_{X_i}+\Delta_i) \cdot R \leq 0$ for an extremal ray $R$ of $\overline{\NE}(X_i/Z)$, 
then $(K_{X_i}+D_i) \cdot R=0$, 
where $\Delta_i$ denotes the push-forward of $\Delta$ on $X_i$. 
\end{enumerate}
\end{prop}

\begin{proof}
Since we have proved the cone theorem in the relative setting (Theorem \ref{t-lc-cone}), 
we can apply the same argument as in \cite[Proposition 3.8]{BW17} (cf.\ \cite[Proposition A.3]{Tanc}). 
\end{proof}

\begin{thm}\label{t-lc-big-bpf2}
Let $k$ be a perfect field of characteristic $p>5$. 
Let $(X, \Delta)$ be a three-dimensional log canonical pair over $k$ and 
let $f:X \to Z$ be a projective $k$-morphism to a quasi-projective $k$-scheme $Z$. 
Let $L$ be an $f$-nef and $f$-big $\R$-Cartier $\R$-divisor such that 
$L-(K_X+\Delta)$ is $f$-semi-ample. 
Then $L$ is $f$-semi-ample. 
\end{thm}

\begin{proof}
We may assume $k$ is algebraically closed.  
By Lemma \ref{l-lcklt-perturb} (1), we may assume that $L=K_X+\Delta$. 
By Corollary \ref{c-dlt-modif}, 
the problem is reduced to the case when $(X, \Delta)$ is a $\Q$-factorial dlt pair. 
Applying Proposition \ref{p-polytope}(3), we can find 
effective $\Q$-divisors $\Delta_1, \ldots, \Delta_n$ and positive real numbers $r_1, \ldots, r_n$ such that 
$(X, \Delta _i)$ is log canonical, $K_X+\Delta_i$ is $f$-nef for any $i \in \{1, \ldots, n \}$, 
\[
\sum_{i=1}^n r_i=1,\quad \text{and}\quad\quad K_X+\Delta=\sum_{i=1}^n r_i(K_X+\Delta_i).
\]
Replacing $\Delta_i$ and $r_i$ appropriately, 
we may assume that each $K_X+\Delta_i$ is $f$-big. 
Hence, each $K_X+\Delta_i$ is $f$-semi-ample by Proposition \ref{p-lc-big-bpf}. 
Therefore, also $L=K_X+\Delta$ is $f$-semi-ample. 
\end{proof}

\begin{thm}\label{t-R-lc-flip}
Let $k$ be a perfect field of characteristic $p>5$. 
Let $(X, \Delta)$ be a three-dimensional log canonical pair over $k$ and 
let $f:X \to Z$ be a projective $k$-morphism to a quasi-projective $k$-scheme $Z$. 
If $K_X+\Delta$ is $f$-big, then there exists a log canonical model of $(X, \Delta)$ over $Z$. 
\end{thm}

\begin{proof}
By Corollary \ref{c-dlt-modif}, 
we may assume that $(X, \Delta)$ is a $\Q$-factorial dlt pair. 
Furthermore, by running a $(K_X+\Delta)$-MMP over $Z$ (Theorem \ref{t-lc-eff-mmp}), 
we may assume that $K_X+\Delta$ is $f$-nef. 
Then Theorem \ref{t-lc-big-bpf2} implies that $K_X+\Delta$ is $f$-semi-ample. 
Let $g:X \to Y$ be the birational morphism $g:X \to Y$ over $Z$ with $g_*\MO_X=\MO_Y$ 
induced by $K_X+\Delta$. 
Then $(Y, \Delta_Y:=g_*\Delta)$ is a log canonical model of $(X, \Delta)$ over $Z$. 
\end{proof}

\section{Base point free theorem}

The purpose of this section is to prove the base point free theorem for log canonical threefolds (Theorem \ref{t-R-rel-bpf2}). 
As a consequence, we obtain the contraction theorem (Theorem \ref{t-contraction}). 
We also establish the minimal model program 
for effective log canonical pairs  (Theorem \ref{t-lc-eff-mmp2}), 
which is a generalisation of Theorem \ref{t-lc-eff-mmp}. 

\begin{lem}\label{l-Q-rel-bpf}
Let $k$ be an algebraically closed field of characteristic $p>5$. 
Let $(X, \Delta)$ be a three-dimensional projective $\Q$-factorial dlt pair over $k$, 
where $\Delta$ is a $\Q$-divisor. 
Let $A$ be an effective big $\Q$-divisor on $X$ such that $(X, \Delta+A)$ is log canonical. 
If $K_X+\Delta+A$ is nef, then $K_X+\Delta+A$ is semi-ample.  
\end{lem}

\begin{proof}
The same argument as in \cite[Step 2-4 of the proof of Theorem 1.2]{Wal} 
works without any changes.  
\end{proof}

\begin{thm}\label{t-R-rel-bpf1}
Let $k$ be a perfect field of characteristic $p>5$. 
Let $(X, \Delta)$ be a three-dimensional log canonical pair over $k$ and  
let $f: X \to Z$ be a projective $k$-morphism to a quasi-projective $k$-scheme $Z$. 
Assume that there exist effective $\R$-divisors $\Delta_1$ and $\Delta_2$ such that $\Delta=\Delta_1+\Delta_2$ and $\Delta_2$ is an $f$-big $\R$-Cartier $\R$-divisor. 
If $K_X+\Delta$ is $f$-nef, then $K_X + \Delta$ is $f$-semi-ample. 
\end{thm}

\begin{proof}
By standard argument, the problem is reduced to the case when 
$k$ is algebraically closed and $f_*\MO_X=\MO_Z$. 
The proof is divided into three steps. 

\setcounter{step}{0}

\begin{step}\label{s1-R-rel-bpf1}
The assertion of Theorem \ref{t-R-rel-bpf1} holds, 
if $Z$ is projective over $k$ and both $\Delta_1$ and $\Delta_2$ are $\Q$-divisors. 
\end{step}

\begin{proof}[Proof of Step \ref{s1-R-rel-bpf1}]
By Theorem \ref{t-lc-cone}, we may assume that $Z=\Spec\,k$. 
Applying Corollary \ref{c-dlt-modif} to $(X, \Delta_1)$, 
there exists a projective birational morphism $g:Y \to X$ such that 
\begin{itemize}
\item 
$Y$ is $\Q$-factorial, and 
\item 
if $B$ is defined by $K_Y+B=g^*(K_X+\Delta_1)$, then 
$B$ is effective and $(Y, \{B\})$ is klt. 
\end{itemize}
In particular, it holds that 
\begin{itemize}
\item 
$(Y, 0)$ is klt, and 
\item 
if we define $\Delta_Y$ by $K_Y+\Delta_Y=g^*(K_X+\Delta)$, then 
$\Delta_Y$ is an effective big $\Q$-divisor such that $(Y, \Delta_Y)$ is log canonical. 
\end{itemize}
Therefore, Lemma \ref{l-Q-rel-bpf} implies that $K_Y+\Delta_Y$ is semi-ample, 
hence so is $K_X+\Delta$. 
This completes the proof of Step \ref{s1-R-rel-bpf1}.
\end{proof}

\begin{step}\label{s2-R-rel-bpf1}
The assertion of Theorem \ref{t-R-rel-bpf1} holds, 
if both $\Delta_1$ and $\Delta_2$ are $\Q$-divisors. 
\end{step}

\begin{proof}[Proof of Step \ref{s2-R-rel-bpf1}]
If $Z=\Spec\,k$, then the assertion follows from Step \ref{s1-R-rel-bpf1}. 
Hence, we may assume that $\dim Z\geq 1$. 
In particular, $\dim X_K \leq 2$, where $X_K$ denotes the generic fibre of $f$. 
Applying Corollary \ref{c-dlt-modif} to $(X, \Delta_1)$, 
we may assume that $X$ is $\Q$-factorial and $(X, 0)$ is klt.

Take an open immersion $Z \hookrightarrow \overline{Z}$ 
to a scheme  $\overline{Z}$ projective over $k$. 
By Proposition \ref{p-lc-cpt3}, there exists 
an open immersion $X \hookrightarrow X^{(1)}$ over $\overline Z$ 
to a $\Q$-factorial threefold projective $X^{(1)}$ over $\overline Z$ 
and an effective $\Q$-divisor $\Delta^{(1)} := \overline{\Delta}$ on $X^{(1)}$ such that 
$(X^{(1)}, \Delta^{(1)})$ is log canonical and $(X^{(1)}, 0)$ is klt.

Since $\dim X_K \leq 2$ and $K_{X^{(1)}}+\Delta^{(1)}$ is pseudo-effective over $\overline Z$, 
it follows from Lemma \ref{l-2dim-nonvani} that 
$K_{X^{(1)}}+\Delta^{(1)} \equiv_f D$ for some effective $\R$-divisor $D$ on $X^{(1)}$. 
By Theorem \ref{t-Wal1.8}, there exists a $(K_{X^{(1)}}+\Delta^{(1)})$-MMP over $\overline Z$ that terminates. 
Let $(X^{(2)}, \Delta^{(2)})$ be the end result. 
Since $K_X+\Delta$ is $f$-nef, 
there exists an open immersion $X \hookrightarrow X^{(2)}$ over $\overline Z$ 
such that $\Delta^{(2)}|_X=\Delta$. 
Since $(X^{(2)}, \Delta^{(2)})$ is log canonical, 
$K_{X^{(2)}}+\Delta^{(2)}$ is nef over $\overline Z$, and 
$\Delta^{(2)}$ is a $\Q$-Cartier $\Q$-divisor which is big over $\overline Z$,  
Step \ref{s1-R-rel-bpf1} implies that $K_{X^{(2)}}+\Delta^{(2)}$ is semi-ample over $\overline Z$. 
Restricting to $X$, it holds that $K_X+\Delta$ is semi-ample over $Z$. 
This completes the proof of Step \ref{s2-R-rel-bpf1}.
\end{proof}

\begin{step}\label{s3-R-rel-bpf1}
The assertion of Theorem \ref{t-R-rel-bpf1} holds without any additional assumptions. 
\end{step}

\begin{proof}[Proof of Step \ref{s3-R-rel-bpf1}]
We may assume that $k$ is an algebraically closed field. 
Applying Corollary \ref{c-dlt-modif} to a log canonical pair $(X, \Delta_1)$, 
the problem is reduced to the case when $X$ is $\Q$-factorial and klt. 
In particular, we may assume that $\Delta_1=0$ and $\Delta=\Delta_2$. 

Then Proposition \ref{p-polytope} implies that 
there exist positive real numbers $r_1, \ldots, r_m$ 
and effective $\Q$-divisors $\Delta_1, \ldots, \Delta_m$ on $X$ such that 
\begin{itemize}
\item $\sum_{i=1}^m r_i=1$,  
\item $(X, \Delta _i)$ is log canonical, $K_X+\Delta_i$ is $f$-nef, and 
\item $K_X+\Delta=\sum_{i=1}^m r_i(K_X+\Delta_i)$. 
\end{itemize}
Replacing $r_i$ and $\Delta_i$ appropriately, we may assume that 
$\Delta_i$ is $f$-big for any $i \in \{1, \ldots, m\}$. 
Hence, Step \ref{s2-R-rel-bpf1} implies that $K_X+\Delta_i$ is $f$-semi-ample 
for any $i \in \{1, \ldots, m\}$. 
Therefore, also $K_X+\Delta$ is $f$-semi-ample. 
This completes the proof of Step \ref{s3-R-rel-bpf1}.
\end{proof}
Step \ref{s3-R-rel-bpf1} completes the proof of Theorem \ref{t-R-rel-bpf1}. 
\end{proof}

\begin{thm}\label{t-R-rel-bpf2}
Let $k$ be a perfect field of characteristic $p>5$. 
Let $(X, \Delta)$ be a three-dimensional log canonical pair over $k$ and 
let $f: X \to Z$ be a projective $k$-morphism to a 
quasi-projective $k$-scheme $Z$. 
Let $L$ be an $f$-nef $\mathbb{R}$-Cartier $\mathbb{R}$-divisor on $X$ such that 
$L - (K_X + \Delta)$ is $f$-semi-ample and $f$-big. Then $L$ is $f$-semi-ample. 
\end{thm}

\begin{proof}
We may assume that $k$ is an algebraically closed field. 
Furthermore, Corollary \ref{c-dlt-modif} reduces the problem to the case when $X$ is $\Q$-factorial. 
Lemma \ref{l-lcklt-perturb} (1) implies that there exists an effective $\R$-Cartier $\R$-divisor $A$ such that 
$A \sim_{f, \R} L - (K_X + \Delta)$ and $(X, \Delta+A)$ is log canonical. 
Since $\Delta+A$ is $f$-big, the assertion follows from Theorem \ref{t-R-rel-bpf1}. 
\end{proof}

\begin{thm}\label{t-contraction}
Let $k$ be a perfect field of characteristic $p>5$. 
Let $(X, \Delta)$ be a three-dimensional log canonical pair over $k$ and 
let $f: X \to Z$ be a projective $k$-morphism to a quasi-projective $k$-scheme $Z$. 
Let $R$ be a $(K_X+\Delta)$-negative extremal ray of $\overline{\NE}(X/Z)$. 
Then, 
\begin{enumerate}
\item[(1)] there exists a projective $Z$-morphism $\varphi_R:X \to Y$ 
such that 
\begin{enumerate}
\item[(1.1)] $Y$ is a normal variety over $k$ projective over $Z$, 
\item[(1.2)] $(\varphi_R)_*\MO_X=\MO_Y$, and 
\item[(1.3)] if $C$ is a curve on $X$ such that $f(C)$ is a point, 
then $\varphi_R(C)$ is a point if and only if $[C] \in R$. 
\end{enumerate}
\end{enumerate}
Moreover, if $\varphi_R:X \to Y$ is a projective $Z$-morphism satisfying the properties (1.1)-(1.3), then the following hold. 
\begin{enumerate}
\item[(2)] 
Fix $\mathbb K \in \{\Q, \R\}$. 
If $L$ is a $\mathbb K$-Cartier $\mathbb K$-divisor on $X$ such that $L \equiv_{\varphi_R} 0$, 
then 
there exists a $\mathbb K$-Cartier $\mathbb K$-divisor $L_Y$ on $Y$ such that 
$L \sim_{f, \mathbb{K}} (\varphi_R)^*L_Y$. 
\item[(3)] 
$\rho(Y/Z)=\rho(X/Z)-1$. 
\item[(4)] 
If $X$ is $\mathbb{Q}$-factorial and if either 
\begin{enumerate}
\item[(4.1)] $\dim X>\dim Y$, or 
\item[(4.2)] $\varphi_R$ is birational and $\dim \Ex(\varphi_R)=2$, 
\end{enumerate}
then $Y$ is $\Q$-factorial. 
\end{enumerate}
\end{thm}

\begin{proof}
Let us show (1). 
By Theorem \ref{t-lc-cone2}, there exists an $f$-ample $\R$-Cartier $\R$-divisor $A$ 
such that $L:=K_X+\Delta+A$ is $f$-nef and $\overline{\NE}(X/Z) \cap L^{\perp} =R.$ 
It follows from Theorem \ref{t-R-rel-bpf2} that 
there is a projective $Z$-morphism $\varphi_R:X \to Y$ which satisfies (1.1)-(1.3). 
Hence, (1) holds. 
%Let $\varphi_R:X \to Y$ be a projective $Z$-morphism satisfying the properties (1.1)-(1.3). 

Let us prove (2). 
We first treat the case when $\overline{\NE}(X/Z) \cap L^{\perp}=R$. 
In this case, $L$ or $-L$ is $f$-nef, hence we may assume that $L$ is $f$-nef. 
By Theorem \ref{t-lc-cone}, 
$L-\epsilon (K_X+\Delta)$ is $f$-ample for some $\epsilon \in \Q_{>0}$. 
Thus, $L$ is $f$-semi-ample by Theorem \ref{t-R-rel-bpf2}. 
In particular, $L$ is $\varphi_R$-semi-ample. 
Since $L \equiv_{\varphi_R} 0$, we can find $L_Y$ as in the statement. 
By Theorem \ref{t-lc-cone2}, 
the general case is reduced to the case when $\overline{\NE}(X/Z) \cap L^{\perp}=R$. 
Hence, (2) holds. 

We now prove (3). 
Fix a curve $\Gamma$ contracted by $\varphi_R$. 
By (2), we have an exact sequence: 
\[
0 \to N^1(Y/Z)_{\Q} \to N^1(X/Z)_{\Q} \xrightarrow{ \cdot \Gamma} \Q \to 0, 
\]
where 
\[
N^1(X/Z)_{\Q}:=\frac{\Pic X \otimes_{\Z} \Q}{\equiv_{Z}}, \quad 
N^1(Y/Z)_{\Q}:=\frac{\Pic Y \otimes_{\Z} \Q}{\equiv_{Z}}.
\]
Then the assertion (3) holds by 
$\rho(X/Z)=\dim_{\Q} N^1(X/Z)_{\Q}$ and $\rho(Y/Z)=\dim_{\Q} N^1(Y/Z)_{\Q}$.

The assertion (4) follows from the same argument as in \cite[Corollary 3.18]{KM98}. 
\end{proof}

\begin{thm}\label{t-lc-eff-mmp2}
Let $k$ be a perfect field of characteristic $p>5$. 
Let $(X, \Delta)$ be a three-dimensional log canonical pair over $k$ and 
let $f:X \to Z$ be a projective $k$-morphism to a quasi-projective $k$-scheme $Z$. 
Assume that $K_X+\Delta \equiv_{f, \R} D$ for some effective $\R$-divisor $D$ on $X$. 
Then we can run an arbitrary $(K_X+\Delta)$-MMP and it terminates. 
\end{thm}

\begin{proof}
By Theorem \ref{t-lc-cone2}, Theorem \ref{t-R-lc-flip} and Theorem \ref{t-contraction}, 
we can run an arbitrary $(K_X+\Delta)$-MMP. 
If $X$ is $\Q$-factorial, then any $(K_X+\Delta)$-MMP terminates by ACC for log canonical thresholds \cite[Theorem 1.10]{Bir16} 
and the assumption $K_X+\Delta \equiv_{f, \R} D \geq 0$ (cf.\ \cite[Lemma 3.2]{Bir07}). 
The general case can be reduced to this case by a standard argument 
(cf.\ \cite[Remark 2.9]{Bir12}). 
\end{proof}

\section{MMP for log canonical pairs}

The purpose of this section is to prove the main theorem of this paper (Theorem \ref{t-lc-mmp}). 
We first check that we may run log minimal model programs with scaling 
under mild conditions (Subsection \ref{ss-exist-scaling}). 
Second, we show that the existence of log minimal models implies  
the termination for certain sequences of flips 
(Subsection \ref{ss1-lc-mmp}). 
Third, we prove the existence of log minimal models 
(Subsection \ref{ss2-lc-mmp}). 
Finally, we establish 
the main theorem of this paper (Subsection \ref{ss3-lc-mmp}).

\subsection{Existence of extremal rays for MMP with scaling}\label{ss-exist-scaling}

In this subsection, we check that we can run log minimal model programs with scaling for any three-dimensional log canonical pairs (Theorem \ref{t-mmp-scaling}).

\begin{lem}\label{l-mmp-scaling}
Let $k$ be a perfect field of characteristic $p>5$. 
Let $(X, \Delta)$ be a three-dimensional geometrically integral log pair over $k$. 
Let $f:X \to Z$ be a projective $k$-morphism to a quasi-projective $k$-scheme $Z$. 
Assume that there exist an $\R$-Cartier $\R$-divisor $C$ on $X$ and 
an effective $\R$-Cartier $\R$-divisor $C'$ on $X_{\overline k}:=X \times_k \overline k$ such that if $\Delta_{\overline k}$ and $C_{\overline k}$ denote the pullbacks of $\Delta$ and $C$ to $X_{\overline k}$ respectively, then 
\begin{enumerate}
\item[(a)] 
$(X_{\overline k}, \Delta_{\overline k}+C')$ is log canonical, 
\item[(b)] 
$K_X+\Delta+C$ is $f$-nef, and 
\item[(c)] 
$C_{\overline k} \sim_{Z, \R} C'$. 
\end{enumerate}
Then $K_{X}+\Delta$ is $f$-nef or there exists a $(K_{X}+\Delta)$-negative extremal ray $R$ of $\overline{\NE}(X/Z)$ such that $(K_{X}+\Delta+\lambda C)\cdot R=0$, where $\lambda$ denotes the $f$-nef threshold of $(K_{X}+\Delta,\,C)$. 
\end{lem}

\begin{proof}
By standard arguments, we may assume that $f_*\MO_X=\MO_Z$. 
We divide the proof into three steps. 

\setcounter{step}{0}

\begin{step}\label{s1-mmp-scaling}
If $k$ is an algebraically closed field, then the assertion of Lemma \ref{l-mmp-scaling} holds. 
\end{step}

\begin{proof}[Proof of Step \ref{s1-mmp-scaling}]
By Theorem \ref{t-lc-cone}, 
the same argument as in \cite[Theorem 4.7.3]{Fuj17} works without any changes. 
\end{proof}

\begin{step}\label{s2-mmp-scaling}
There exists a curve $G$ on $X$ such that $f(G)$ is a point, $(K_{X}+\Delta)\cdot G <0$ and $(K_{X}+\Delta+\lambda C)\cdot G=0$. 
\end{step}

\begin{proof}[Proof of Step \ref{s2-mmp-scaling}]
Let $f_{\overline k}:X_{\overline k}\to Z_{\overline k}$ be the base change $f \times_k \overline k$. 
Applying Step \ref{s1-mmp-scaling} to $f_{\overline k}:X_{\overline k}\to Z_{\overline k}$, we can find a curve $G'$ on $X_{\overline k}$ such that $f_{\overline k}(G')$ is a point, $(K_{X_{\overline k}}+\Delta_{\overline k})\cdot G'<0$ and $(K_{X_{\overline k}}+\Delta_{\overline k}+\lambda C_{\overline k})\cdot G'=0$. 
Note that the $f_{\overline k}$-nef threshold of $(K_{X_{\overline k}}+\Delta_{\overline k},\,C_{\overline k})$ is $\lambda$ (cf.~\cite[Remark 2.7]{GNT}).
Let $G$ be the image of $G'$ on $X$. 
Then $f(G)$ is a point,  $(K_{X}+\Delta)\cdot G <0$ and $(K_{X}+\Delta+\lambda C)\cdot G=0$ (cf.~\cite[Lemma 2.3]{Tana}). 
Hence $G$ is the required curve. 
This completes the proof of Step \ref{s2-mmp-scaling}. 
\end{proof}

\begin{step}\label{s3-mmp-scaling}
The assertion of Lemma \ref{l-mmp-scaling} holds without any additional assumptions. 
\end{step}

\begin{proof}[Proof of Step \ref{s3-mmp-scaling}]
Let $G$ be a curve as in Step \ref{s2-mmp-scaling}. 
Let $A$ be an $f$-ample $\R$-Cartier $\R$-divisor on $X$ such that $(K_{X}+\Delta+A)\cdot G <0$. 
By Theorem \ref{t-lc-cone2}, there exist $G_{0}\in \overline{\NE}(X/Z)_{K_{X}+\Delta+A\geq 0}$, 
$r_{1}, \ldots ,r_m\in \R_{\geq 0}$, and 
curves $G_{1}, \ldots ,G_m$ 
generating $(K_{X}+\Delta+A)$-negative extremal rays of $\overline{\NE}(X/Z)$ 
such that the equation 
\[
[G]=[G_{0}]+\sum_{i=1}^m r_{i} [G_{i}]
\]
holds in $\overline{\NE}(X/Z)$. 
Since $(K_{X}+\Delta+A)\cdot G <0$, there is $j \in \{1, \ldots ,\,m\}$ such that $r_j >0$ and $(K_{X}+\Delta+A)\cdot G_{j} <0$. 
In particular, we get $(K_{X}+\Delta)\cdot G_{j} <0$. 
On the other hand, since $K_{X}+\Delta+\lambda C$ is $f$-nef, 
the equation $(K_{X}+\Delta+\lambda C)\cdot G=0$ implies that 
$(K_{X}+\Delta+\lambda C)\cdot G_{j}=0$. 
Therefore $G_{j}$ generates a $(K_{X}+\Delta)$-negative extremal ray $R$ 
of $\overline{\NE}(X/Z)$ such that $(K_{X}+\Delta+\lambda C)\cdot R=0$. 
This completes the proof of Step \ref{s3-mmp-scaling}. 
\end{proof}
Step \ref{s3-mmp-scaling} completes the proof of Lemma \ref{l-mmp-scaling}. 
\end{proof}

By Lemma \ref{l-mmp-scaling}, we can run an MMP under the same assumption. 
However, for now, we do not know whether it terminates.

\begin{thm}\label{t-mmp-scaling}
Let $k$ be a perfect field of characteristic $p>5$. 
Let $(X, \Delta)$ be a three-dimensional geometrically integral log pair over $k$. 
Let $f:X \to Z$ be a projective $k$-morphism to a quasi-projective $k$-scheme $Z$. 
Assume that there exist an $\R$-Cartier $\R$-divisor $C$ on $X$ and 
an effective $\R$-Cartier $\R$-divisor $C'$ on $X_{\overline k}:=X \times_k \overline k$ such that if $\Delta_{\overline k}$ and $C_{\overline k}$ denote the pullbacks of $\Delta$ and $C$ to $X_{\overline k}$ respectively, then 
\begin{enumerate}
\item[(a)] 
$(X_{\overline k}, \Delta_{\overline k}+C')$ is log canonical, 
\item[(b)] 
$K_X+\Delta+C$ is $f$-nef, and 
\item[(c)] 
$C_{\overline k} \sim_{Z, \R} C'$. 
\end{enumerate}
Then there exists a $(K_{X}+\Delta)$-MMP over $Z$ with scaling of $C$ 
\begin{equation}\label{e1-mmp-scaling}
X=:X_0 \dashrightarrow X_1 \dashrightarrow \cdots 
\end{equation}
such that either 
\begin{enumerate}
\item (\ref{e1-mmp-scaling}) terminates, or 
\item (\ref{e1-mmp-scaling}) is an infinite sequence. 
\end{enumerate}
\end{thm}

\begin{proof}
The assertion follows from Theorem \ref{t-R-lc-flip}, Theorem \ref{t-contraction}, 
and Lemma \ref{l-mmp-scaling} (cf.\ \cite[4.9.1]{Fuj17}). 
\end{proof}

\begin{rem}
Let $k$, $(X, \Delta)$ and $f:X \to Z$ be as in Theorem \ref{t-mmp-scaling}. 
If $C$ is a sufficiently large multiple of an $f$-ample $\R$-Cartier $\R$-divisor on $X$, 
then there exists $C'$ on $X_{\overline k}$ that satisfies all the conditions 
(a), (b) and (c) of Theorem \ref{t-mmp-scaling}. 
\end{rem}

\subsection{Criterion for termination of flips}\label{ss1-lc-mmp}

In this subsection, we prove that assuming the existence of log minimal models, 
the termination holds for pseudo-effective minimal model programs with scaling 
whose scaling coefficients are strictly decreasing. 
The idea of the proof can be found in the proof of \cite[Theorem 4.1(iii)]{Bir12}.

\begin{prop}\label{p-dlt-mmp3}
Let $k$ be an algebraically closed field of characteristic $p>5$. 
Let $X$ be a $\Q$-factorial normal threefold over $k$ and 
let $f:X \to Z$ be a projective $k$-morphism to a quasi-projective $k$-scheme $Z$. 
Let $\Delta$ and $C$ be effective $\R$-divisors on $X$ such that 
\begin{enumerate}
\item[(a)] 
$(X, \Delta+C)$ is log canonical, 
\item[(b)] 
$K_X+\Delta+C$ is $f$-nef, and 
\item[(c)] 
$C$ is $f$-big. 
\end{enumerate}
Assume that the following holds: 
\begin{enumerate}
\item[(i)]  
If $(V, \Delta_V)$ is a three-dimensional $\Q$-factorial log canonical pair 
that is projective over $Z$ and $K_V+\Delta_V$ is pseudo-effecitve over $Z$, 
then there exists a log minimal model of $(V, \Delta_V)$ over $Z$. 
\end{enumerate}
Then there exists no infinite sequence 
\begin{equation}\label{e1-dlt-mmp3}
X=X_0 \dashrightarrow \cdots \dashrightarrow X_i \dashrightarrow X_{i+1}\dashrightarrow \cdots
\end{equation}
such that 
\begin{enumerate}
\item[(ii)] 
the sequence (\ref{e1-dlt-mmp3}) is a $(K_X+\Delta)$-MMP over $Z$ with scaling of $C$, and  
\item[(iii)] 
if $\lambda_0, \lambda_1, \ldots$ are the real numbers defined by 
\[
\lambda_i:=\min\{\mu \in \R_{\geq 0} \,|\, K_{X_i}+\Delta_i+\mu C_i
{\rm \,\,is\,\,nef}\},
\] 
then it holds that $\lim_{i \to \infty}\lambda_i\neq \lambda_i$ for any $i$. 
\end{enumerate}
\end{prop}

\begin{proof}
Assume that there exists an infinite sequence (\ref{e1-dlt-mmp3}) which satisfies (ii) and (iii). 
Let us derive a contradiction. 

Set $\lambda_{\infty}:=\lim_{i \to \infty} \lambda_i$. 
Replacing $\Delta$ and $C$ 
by $\Delta+\lambda_{\infty} C$ and $(1-\lambda_{\infty})C$ respectively, 
we may assume that 
\begin{enumerate}
\item[(iii)'] $\lambda_{\infty}=\lim_{i \to \infty} \lambda_i=0$, and 
\end{enumerate}
\begin{enumerate}
\item[(iv)\;] 
$K_X+\Delta$ is $f$-pseudo-effective. 
\end{enumerate}

It is clear that $\lambda_{0}>0$. 
By (i), there is a log minimal model $(X',\Delta')$ of $(X,\Delta)$ over $Z$. 
Let $g:Y\to X$ be a log resolution of $(X,\Delta+C)$ such that the induced birational map $h:Y \dashrightarrow X'$ is a morphism. 
Set $\Gamma:=g_*^{-1}\Delta+E_g$, 
where $E_g$ denotes the reduced divisor on $Y$ such that $\Supp\,E_g=\Ex(g)$. 
By construction, $(Y,\Gamma)$ is a log birational model over $Z$ of $(X,\Delta)$ (cf.~Definition \ref{d-lmm}) and $(Y, \Gamma+g_{*}^{-1}C)$ is log canonical. 
Moreover, by a property of weak log canonical model, we can write 
$$K_{Y}+\Gamma=h^{*}(K_{X'}+\Delta')+E$$
with an $h$-exceptional divisor $E\geq0$. 
It follows from Theorem \ref{t-lc-eff-mmp} that 
there exists a $(K_{Y}+\Gamma)$-MMP over $X'$ that terminates 
with a log minimal model $(Y', \Gamma')$ over $X'$. 
The negativity lemma implies that $K_{Y'}+\Gamma'=h'^{*}(K_{X'}+\Delta')$, 
where $h':Y' \to X'$ denotes the induced birational morphism. 
In particular $K_{Y'}+\Gamma'$ is nef over $Z$. 

Pick a sufficiently small $\epsilon>0$ so that 
the MMP $Y \dashrightarrow Y'$ defined above is 
a $(K_{Y}+\Gamma+\epsilon g_{*}^{-1}C)$-MMP over $Z$. 
Let $C'$ be the proper transform of $g_{*}^{-1}C$ on $Y'$. 
Then $(Y',\Gamma'+\epsilon C')$ is a $\mathbb{Q}$-factorial log canonical pair. 
Since it follows from (iv) that $K_Y+\Gamma+\epsilon g_{*}^{-1}C$ is big over $Z$, 
so is $K_{Y'}+\Gamma'+\epsilon C'$. 
Thus, again by Theorem \ref{t-lc-eff-mmp}, 
there is a $(K_{Y'}+\Gamma'+\epsilon C')$-MMP over $Z$ that terminates 
with a log minimal model $(Y'',\Gamma''+\epsilon C'')$ over $Z$.  
As $\epsilon$ is sufficiently small, 
it follows from Proposition \ref{p-polytope} that 
the MMP $Y' \dashrightarrow Y''$ is $(K_{Y'}+\Gamma')$-numerically trivial. 
Therefore the $\R$-divisor $K_{Y''}+\Gamma''$ is also nef over $Z$. 
Hence, we see that $K_{Y''}+\Gamma''+\epsilon'C''$ is nef over $Z$ 
for any $\epsilon' \in [0,\epsilon]$. 

By (iii)', there exists $i>0$ such that $0\leq \lambda_{i+1}<\lambda_{i}\leq \epsilon$. 
By the construction of $(Y,\Gamma)$ and by the basic property of the log MMP, 
we see that $(X_{i+1},\Delta_{i+1}+\lambda_{i}C_{i+1})$ and  $(X_{i+1},\Delta_{i+1}+\lambda_{i+1}C_{i+1})$ are weak log canonical models over $Z$ of 
$(Y,\Gamma+\lambda_{i}g_{*}^{-1}C)$ and $(Y,\Gamma+\lambda_{i+1}g_{*}^{-1}C)$ respectively. 
On the other hand, by the construction of $(Y'',\Gamma'')$ and by the choices of $\epsilon$ and $i$, 
we see that also $(Y'',\Gamma''+\lambda_{i}C'')$ and $(Y'',\Gamma''+\lambda_{i+1}C'')$ are weak log canonical models over $Z$ of 
$(Y,\Gamma+\lambda_{i}g_{*}^{-1}C)$ and $(Y,\Gamma+\lambda_{i+1}g_{*}^{-1}C)$ respectively. 
Now let $\varphi:W\to X_{i+1}$ and $\psi:W\to Y''$ be a common resolution of $X_{i+1}\dashrightarrow Y''$. 
Then it follows from Remark \ref{r-2-wlc-models} that 
\begin{equation*}
\begin{split}
\varphi^{*}(K_{X_{i+1}}+\Delta_{i+1}+\lambda_{i}C_{i+1})&=\psi^{*}(K_{Y''}+\Gamma''+\lambda_{i}C'')\qquad{\rm and }\\
\varphi^{*}(K_{X_{i+1}}+\Delta_{i+1}+\lambda_{i+1}C_{i+1})&=\psi^{*}(K_{Y''}+\Gamma''+\lambda_{i+1}C'').
\end{split}
\end{equation*}
By $\lambda_{i}\neq \lambda_{i+1}$, we get 
$$\varphi^{*}(K_{X_{i+1}}+\Delta_{i+1})=\psi^{*}(K_{Y''}+\Gamma'').$$
Therefore, $K_{X_{i+1}}+\Delta_{i+1}$ is nef over $Z$. 
Hence, the $(K_{X}+\Delta)$-MMP (\ref{e1-dlt-mmp3}) terminates. 
This contradicts the fact that the sequence (\ref{e1-dlt-mmp3}) 
was chosen to be an infinite sequence. 
\end{proof}

\subsection{Existence of log minimal models}\label{ss2-lc-mmp}

In this subsection, we prove the existence of log minimal models (Theorem \ref{t-rel-lmm}). 
To this end, we first treat the projective case (Theorem \ref{t-exist-lmm}). 
The proof of Theorem \ref{t-exist-lmm} is similar to the one of \cite[Corollary 1.7]{Bir12b}. 
We start with an auxiliary result, which is known to experts.

\begin{lem}\label{l-terminal-flip}
Let $k$ be an algebraically closed field of characteristic $p>5$. 
Let $(X, \Delta)$ be a three-dimensional $\Q$-factorial terminal pair over $k$ and 
let $f:X \to Z$ be a projective $k$-morphism to a quasi-projective $k$-scheme $Z$. 
Then there exists no infinite sequence 
\begin{equation}\label{e1-terminal-flip}
X=X_0 \overset{\varphi_0}{\dashrightarrow} X_1 \overset{\varphi_1}{\dashrightarrow} \cdots
\end{equation}
such that 
\begin{enumerate}
\item the sequence $(\ref{e1-terminal-flip})$ is a $(K_X+\Delta)$-MMP over $Z$, and 
\item for any $i$, $\varphi_i:X_i \dashrightarrow X_{i+1}$ is a $(K_{X_i}+\Delta_i)$-flip, 
where $\Delta_i$ denotes the proper transform of $\Delta$ to $X_i$. 
\end{enumerate}
\end{lem}

\begin{proof}
By the fact that the singular locus of $X$ is zero-dimensional \cite[Corollary 2.30]{Kol13}, 
we can apply the same argument as in \cite[Theorem 6.17]{KM98}. 
\end{proof}

\begin{thm}\label{t-exist-lmm}
Let $k$ be an algebraically closed field of characteristic $p>5$. 
Let $(X, \Delta)$ be a three-dimensional projective log canonical pair over $k$ such that $K_X+\Delta$ is pseudo-effective. 
Then there exists a log minimal model of $(X, \Delta)$.  
\end{thm}

\begin{proof}
Set $S:=\lfloor \Delta \rfloor$ and $B:=\{\Delta\}$. 
In particular, we have $\Delta=S+B$. 
Taking a log resolution of $(X, \Delta)$, we may assume that 
\begin{enumerate}
\item 
$(X, \Delta)$ is a $\Q$-factorial dlt pair. 
\end{enumerate}
If there exists a $(K_X+\Delta)$-MMP that terminates, then there is nothing to prove. 
Hence, the problem is reduced to the case when 
\begin{enumerate}
\item[(2)]  
an arbitrary $(K_X+\Delta)$-MMP does not terminate. 
\end{enumerate}
We divide the proof into several steps. 

\setcounter{step}{0}

\begin{step}\label{s1-exist-lmm}
There exist an effective big $\R$-divisor $H$ on $X$ and an infinite sequence 
\begin{equation}\label{e1-exist-lmm}
X=X_0 \overset{f_0}{\dashrightarrow} X_1\overset{f_1}{\dashrightarrow}  \cdots 
\end{equation}
such that 
\begin{enumerate}
\item[(3)] $K_X+\Delta+H$ is nef, $(X, \Delta + H)$ is dlt,  
\item[(4)] 
the sequence (\ref{e1-exist-lmm}) is a $(K_X+\Delta)$-MMP with scaling of $H$, and  
\item[(5)] 
$\lim_{n \to \infty} \lambda_n =0$, where $\lambda_n$ is the scaling coefficient. 
\end{enumerate}
\end{step}

\begin{proof}[Proof of Step \ref{s1-exist-lmm}]
Pick an ample $\R$-divisor $H$ such that $K_{X}+\Delta+H$ is nef and 
$(X, \Delta + H)$ is dlt \cite[Lemma 9.2]{Bir16}. 
Hence, (3) holds. 
By Theorem \ref{t-mmp-scaling} and (2), 
there exists an infinite sequence which is a $(K_X+\Delta)$-MMP with scaling of $H$: 
\begin{equation}\label{e2-exist-lmm2}
(X, \Delta)=(X_{0}, \Delta_{0})\overset{f_0}{\dashrightarrow} (X_{1},\Delta_{1})\overset{f_1}{\dashrightarrow} \cdots \dashrightarrow (X_{i},\Delta_{i})\dashrightarrow \cdots .
\end{equation}
Clearly, (4) holds. 

It suffices to prove (5). 
Assuming that $\lambda:=\lim_{n \to \infty} \lambda_n > 0$, let us derive a contradiction. 
The sequence (\ref{e2-exist-lmm2}) is a $(K_X+\Delta+\frac{\lambda}{2}H)$-MMP, 
thus there exists a $(K_X+\Delta+\frac{\lambda}{2}H)$-MMP that does not terminate. 
On the other hand, $K_X+\Delta+\frac{\lambda}{2}H$ is big, 
hence any $(K_X+\Delta+\frac{\lambda}{2}H)$-MMP terminates by Theorem \ref{t-lc-eff-mmp2}. 
This is a contradiction. 
Therefore, (5) holds. 
\end{proof}

We denote the proper transforms of $\Delta, S, B$ and $H$ on $X_i$ by 
$\Delta_i, S_i, B_i$ and $H_i$, respectively. 
Note that even if we replace $(X, \Delta)$ by $(X_i, \Delta_i)$, 
the properties (1)--(5) still hold. 
In particular, we may assume that 

\begin{enumerate}
\item[(6)] $1 > \lambda_0 \geq \lambda_1 \geq \cdots$, 
\item[(7)] for any $i$, $X_i \dashrightarrow X_{i+1}$ is a flip, and 
\item[(8)] for any $i$, $\Ex(f_i)$ is disjoint from $\Supp\,S_i$  (cf.\ \cite[Proposition 5.5]{Bir16}, \cite[Proposition 4.2]{Wal}). 
\end{enumerate}
Hence, $(X_{i},B_{i}+\lambda_{i}H_{i})$ is klt, because $(X,B+\lambda_{i}H)$ is klt 
and the MMP 
\[
X=X_0 \dashrightarrow X_1 \dashrightarrow \cdots \dashrightarrow X_i
\]
is $(K_X+B+\lambda_i H)$-non-positive. 
For any $i\geq0$, let $\mu_{i}:(W_{i},\Psi_{i}) \to (X_{i},B_{i}+\lambda_{i}H_{i})$ be a projective birational morphism such that $(W_{i},\Psi_{i})$ is a $\Q$-factorial terminal pair and $K_{W_{i}}+\Psi_{i}=\mu_{i}^{*}(K_{X_{i}}+B_{i}+\lambda_{i}H_{i})$. 
Let $g_i:W_i \dashrightarrow W_{i+1}$ be the induced birational map. 
Since $a_D(X_{i},B_{i}+\lambda_{i}H_{i}) \leq a_D(X_{i+1},B_{i+1}+\lambda_{i+1}H_{i+1})$ 
for any exceptional prime divisor $D$ over $X$, 
the induced birational map $g_i^{-1}:W_{i+1} \dashrightarrow W_i$ 
does not contract any prime divisor. 
Replacing  $(X,\Delta)$ by $(X_{i},\Delta_{i})$ for some $i \gg 0$, 
we may assume that 
\begin{enumerate}
\item[(9)]  
$g_i:W_{i}\dashrightarrow W_{i+1}$ is isomorphic in codimension one for any $i$. 
\end{enumerate} 
Set $h_{i}:W_{0}\dashrightarrow W_{i}$ to be the induced birational map. 
Then we get 
$$\Psi_{0}\geq h_{1*}^{-1}\Psi_{1}\geq\cdots \geq h_{i*}^{-1}\Psi_{i} \geq \cdots \geq 0.$$
Therefore there exists an $\R$-divisor $\Psi_{\infty}$ on $W_0$ such that 
$\Psi_{\infty}={\rm lim}_{i\to \infty}h_{i*}^{-1}\Psi_{i}$. 
Then $(W_{0},\Psi_{\infty})$ is a $\mathbb{Q}$-factorial terminal pair. 
Set 
\[
G_i=\mu_{i}^{*}S_i.
\] 

\begin{step}\label{s2-exist-lmm}
The following hold. 
\begin{enumerate}
\item[(i)] $(W_{0},G_0 +\Psi_{\infty})$ is log canonical. 
\item[(ii)] If there exists a log minimal model of $(W_0, G_0+\Psi_{\infty})$, 
then there exists a log minimal model of $(X, \Delta)$. 
\end{enumerate}
\end{step}

\begin{proof}[Proof of Step \ref{s2-exist-lmm}]
If we define $\Xi_{i}$ on $W_{0}$ by $K_{W_{0}}+\Xi_{i}=\mu_{0}^{*}(K_{X}+B+\lambda_{i}H)$, then 
$\Xi_{i}\geq h _{i*}^{-1}\Psi_{i}$ by the negativity lemma. 
Since ${\rm lim}_{i\to \infty}\lambda_{i}=0$, we have 
$K_{W_{0}}+G_0 +\Psi_{\infty}\leq \mu _{0}^{*}(K_{X}+\Delta).$ 
Therefore $(W_{0},G_0 +\Psi_{\infty})$ is log canonical. 

If $(W_{0},G_0 +\Psi_{\infty})$ has a log minimal model, 
then $K_{X}+\Delta$ has a weak Zariski decomposition 
in the sense of \cite[Section 8.1]{Bir16}. 
So we see that $(X,\Delta)$ has a log minimal model by \cite[Proposition 8.3]{Bir16}. 
\end{proof}

\begin{step}\label{s3-exist-lmm}
The following hold. 
\begin{enumerate}
\item[(i)] 
$G_0=(h_i)^{-1}_*G_i$ for any $i$. 
\item[(ii)] 
$K_{W_0}+G_0+\Psi_{\infty}=\lim_{i \to \infty} (h_i)_*^{-1}(K_{W_i}+G_i+\Psi_i).$
\item[(iii)] 
$K_{W_0}+G_0+\Psi_{\infty} \in \overline{\Mov}(X)$. 
\item[(iv)] 
The stable base locus of $(h_i)_*^{-1}(K_{W_i}+G_i+\Psi_i)$ is disjoint from $\Supp\,G_0$. 
\end{enumerate}
\end{step}

\begin{proof}[Proof of Step \ref{s3-exist-lmm}]
The assertion (i) follows from (8) and (9). 
Then (ii) holds by (i). 

We now show (iii). 
Since $K_{X_i}+\Delta_i+\lambda_i H_i$ is nef and big, 
it follows from Theorem \ref{t-Wal1.1} that $K_{X_i}+\Delta_i+\lambda_i H_i$ is semi-ample. 
Since 
\[
K_{W_i}+G_i+\Psi_i=\mu_i^* (K_{X_i}+S_i+B_i+\lambda_i H_i)=\mu_i^* (K_{X_i}+\Delta_i+\lambda_i H_i), 
\]
also $K_{W_i}+G_i+\Psi_i$ is semi-ample. 
Then (ii) implies (iii).  

Let us prove (iv). 
Let $\varphi_{0}:Y \to W_{0}$ and $\varphi_{i}:Y \to W_{i}$ be a common resolution of 
$W_{0}$ and $W_{i}$. 
Since $K_{W_{i}}+G_{i}+\Psi_{i}$ is nef, 
the negativity lemma induces an equation 
\[
\varphi_{0}^{*}((h_i)_*^{-1}(K_{W_{i}}+G_{i}+\Psi_{i}))=\varphi_{i}^{*}(K_{W_{i}}+G_{i}+\Psi_{i})+F
\]
for some $\varphi_{0}$-exceptional effective $\R$-divisor $F$. 
Since $K_{W_{i}}+G_{i}+\Psi_{i}$ is semi-ample, the stable base locus of $(h_i)_*^{-1}(K_{W_{i}}+G_{i}+\Psi_{i})$ is $\varphi_{0}({\rm Supp}\,F)$. 
By (8), there is an open set $U_{i}\subset X$ containing $S$ such that the induced birational map $X\dashrightarrow X_{i}$ is an isomorphism on $U_{i}$. 
Restricting the above equation to $(\mu_0 \circ \varphi_0)^{-1}(U_{i})$, 
the negativity lemma implies that $F\!\mid\!_{(\mu_0 \circ \varphi_0)^{-1}(U_{i})}=0$. 
Thus we have $\varphi_{0}({\rm Supp}\,F)\cap \mu_{0}^{-1}(U_{i})=\emptyset$. 
Since $\mu_{0}^{-1}(U_{i})$ contains ${\rm Supp}\,G_0$, we see that the stable base locus of 
$(h _{i})_{*}^{-1}(K_{W_{i}}+G_{i}+\Psi_{i})$, which is equal to $\varphi_0(\Supp\,F)$, is disjoint from ${\rm Supp}\, G_0$. 
\end{proof}

\begin{step}\label{s4-exist-lmm}
There exists a log minimal model of $(W_0, G_0+\Psi_{\infty})$. 
\end{step}

\begin{proof}[Proof of Step \ref{s4-exist-lmm}]
Let 
\begin{equation}\label{e2-exist-lmm}
W_0=:V_0 \dashrightarrow V_1 \dashrightarrow \cdots 
\end{equation}
be an arbitrary $(K_{W_0}+G_0+\Psi_{\infty})$-MMP. 
It suffices to prove that the MMP (\ref{e2-exist-lmm}) terminates. 
By Step \ref{s3-exist-lmm}(iii), 
the MMP (\ref{e2-exist-lmm}) consists of flips. 
Furthermore, Step \ref{s3-exist-lmm}(ii)(iv) imply that the MMP (\ref{e2-exist-lmm}) 
occurs disjointly from $\Supp\,G_0$. 
In particular, the sequence (\ref{e2-exist-lmm}) is a $(K_{W_0}+\Psi_{\infty})$-MMP. 
Since $(W_0, \Psi_{\infty})$ is a $\Q$-factorial terminal pair, 
the sequence (\ref{e2-exist-lmm}) terminates by Lemma \ref{l-terminal-flip}. 
\end{proof}
Step \ref{s2-exist-lmm} and Step \ref{s4-exist-lmm} complete the proof of Theorem \ref{t-exist-lmm}. 
\end{proof}

\begin{prop}\label{p-abs-mmp}
Let $k$ be an algebraically closed field of characteristic $p>5$. 
Let $(X, \Delta)$ be a three-dimensional $\Q$-factorial log canonical pair over $k$ 
such that $(X, \{\Delta\})$ is klt. 
Let $f:X \to Z$ be a projective $k$-morphism to a projective $k$-scheme $Z$. 
Let $C$ be an $f$-ample $\R$-Cartier $\R$-divisor such that $K_{X}+\Delta+C$ is $f$-nef. 
Then the following hold. 
\begin{enumerate}
\item 
There exists no infinite sequence that is a $(K_X+\Delta)$-MMP over $Z$ with scaling of $C$. 
\item 
There exists a $(K_X+\Delta)$-MMP over $Z$ with scaling of $C$ that terminates. 
\end{enumerate}
\end{prop}

\begin{proof}
By Theorem \ref{t-mmp-scaling}, 
(1) implies (2).

Let us prove (1). 
Assume that there exists an infinite sequence
\begin{equation}\label{e1-abs-mmp}
X=X_0 \dashrightarrow X_1 \dashrightarrow \cdots 
\end{equation}
 that is a $(K_X+\Delta)$-MMP over $Z$ with scaling of $C$. 
Let us derive a contradiction. 
If the scaling coefficients satisfy $\epsilon := \lim \lambda_i>0$, 
then the sequence (\ref{e1-abs-mmp}) is a $(K_X+\Delta+\epsilon C)$-MMP 
with scaling of $C$. 
Since there exists an effective $\R$-divisor $\Delta'$ such that 
$\Delta' \sim_{Z, \R} \Delta+\epsilon C$ and $(X, \Delta')$ is klt, 
it follows from \cite[Theorem 1.5]{BW17} 
that the sequence (\ref{e1-abs-mmp}) terminates, which is a contradiction. 
Hence, we may assume that $\lim \lambda_i =0$. 
By Theorem \ref{t-exist-lmm}, 
the sequence (\ref{e1-abs-mmp}) terminates by Proposition \ref{p-dlt-mmp3}. 
Therefore, we get a contraction in any case. Hence, (1) holds. 
\end{proof}

\begin{cor}\label{c-abs-mmp}
Let $k$ be an algebraically closed field of characteristic $p>5$. 
Let $(X, \Delta)$ be a three-dimensional $\Q$-factorial log canonical pair over $k$ 
such that $(X, \{\Delta\})$ is klt. 
Let $f:X \to Z$ be a projective $k$-morphism to a quasi-projective $k$-scheme $Z$. 
Then there exists a $(K_X+\Delta)$-MMP over $Z$ that terminates. 
\end{cor}

\begin{proof}
The assertion follows from 
Proposition \ref{p-lc-cpt} and Proposition \ref{p-abs-mmp}. 
\end{proof}

\begin{thm}\label{t-rel-lmm}
Let $k$ be an algebraically closed field of characteristic $p>5$. 
Let $(X, \Delta)$ be a three-dimensional log canonical pair over $k$ and 
let $f:X \to Z$ be a projective $k$-morphism to a quasi-projective $k$-scheme $Z$. 
If $K_X+\Delta$ is $f$-pseudo-effective, then there exists a 
log minimal model of $(X, \Delta)$ over $Z$. 
\end{thm}

\begin{proof}
If $(X, \Delta)$ is $\Q$-factorial dlt pair, then the assertion follows from 
Corollary \ref{c-abs-mmp}. 
The general case is reduced to this case by Corollary \ref{c-dlt-modif}. 
\end{proof}

\begin{cor}\label{c-rel-lmm}
Let $k$ be a perfect field of characteristic $p>5$. 
Let $(X, \Delta)$ be a three-dimensional geometrically integral log pair over $k$. 
Let $f:X \to Z$ be a projective $k$-morphism to a quasi-projective $k$-scheme $Z$. 
Assume that there exist an $\R$-Cartier $\R$-divisor $C$ on $X$ and 
an effective $\R$-Cartier $\R$-divisor $C'$ on $X_{\overline k}:=X \times_k \overline k$ such that if $\Delta_{\overline k}$ and $C_{\overline k}$ denote the pullbacks of $\Delta$ and $C$ to $X_{\overline k}$ respectively, then 
\begin{enumerate}
\item[(a)] 
$(X_{\overline k}, \Delta_{\overline k}+C')$ is log canonical, 
\item[(b)] 
$K_X+\Delta+C$ is $f$-nef, 
\item[(c)] 
$C$ is $f$-big, and  
\item[(d)] 
$C_{\overline k} \sim_{Z, \R} C'$. 
\end{enumerate}
Then there exists no infinite sequence 
\begin{equation}\label{e1-dlt-mmp2}
X=X_0 \dashrightarrow \cdots \dashrightarrow X_i \dashrightarrow X_{i+1}\dashrightarrow \cdots
\end{equation}
such that 
\begin{enumerate}
\item[(1)] 
the sequence (\ref{e1-dlt-mmp2}) is a $(K_X+\Delta)$-MMP over $Z$ with scaling of $C$, and  
\item[(2)] 
if $\lambda_0, \lambda_1, \ldots$ are the real numbers defined by 
\[
\lambda_i:=\min\{\mu \in \R_{\geq 0} \,|\, K_{X_i}+\Delta_i+\mu C_i
{\rm \,\,is\,\,nef\,\,over\,\,}Z\},
\] 
then it holds that $\lim_{i \to \infty}\lambda_i\neq \lambda_i$ for any $i$. 
\end{enumerate}
\end{cor}

\begin{proof}
Under the additional assumption that 
\begin{itemize}
\item $k$ is algebraically closed and $X$ is $\mathbb{Q}$-factorial and klt, 
\end{itemize}
the assertion immediately follows from Proposition \ref{p-dlt-mmp3} 
and Theorem \ref{t-rel-lmm}. 

Let us go back to the general situation. 
Let 
\begin{equation*}
X'=X'_0 \dashrightarrow X'_1 \dashrightarrow \cdots 
\end{equation*}
be the sequence obtained by the base changes $X'_i:=X_i \times_k \overline k$ 
to the algebraic closure $\overline k$ of $k$. 
Then each birational map $g'_i:X'_i \dashrightarrow X'_{i+1}$ has a decomposition 
\[
X'_i \xrightarrow{\varphi_i} Z'_i \xleftarrow{\psi_i} X'_{i+1}, 
\] 
such that 
\begin{enumerate}
\item[(i)] $Z'_i$ is a normal threefold projective over $Z$, 
\item[(ii)] $\varphi_i$ and $\psi_i$ are birational, 
\item[(iii)] $-(K_{X'_i}+\Delta'_i)$ is $\varphi_i$-ample, and 
\item[(iv)] $(X'_{i+1}, \Delta'_{i+1})$ is a log canonical model of $(X'_i, \Delta'_i)$ over $Z'_i$.  
\end{enumerate}
We apply Corollary \ref{c-dlt-modif} to $(X'_0, \Delta'_0)$. 
Then we obtain a projective birational morphism $g_0:Y_0 \to X'_0$ 
satisfying the properties listed in Corollary \ref{c-dlt-modif}. 
In particular, if we define $\Delta_{Y_0}$ by $K_{Y_0}+\Delta_{Y_0}=g_0^*(K_{X'_0}+\Delta'_0)$, then $(Y_{0}, \Delta_{Y_0})$ is a $\Q$-factorial dlt pair. 
Then there exists a $(K_{Y_0}+\Delta_{Y_0})$-MMP over $Z'_0$ that terminates. 
Let $(Y_1, \Delta_{Y_1})$ be the end result, which is a log minimal model of $(X'_0, \Delta'_0)$ over $Z'_0$. 
By (iv) and Remark \ref{r-2-wlc-models}, 
the induced birational map $g_1:Y_1 \dashrightarrow X'_1$ is a morphism and $K_{Y_1}+\Delta_{Y_1}=g_1^*(K_{X'_1}+\Delta'_1)$. 
Repeating the same procedure, we obtain a sequence 
\[
Y_0 \dashrightarrow Y_1 \dashrightarrow Y_2 \dashrightarrow \cdots, 
\]
which is a $(K_{Y_0}+\Delta_{Y_0})$-MMP over $Z' := Z \times _k \overline{k}$. 
Furthermore, the MMP $Y_i \dashrightarrow Y_{i+1}$ is $(K_{Y_i}+\Delta_{Y_i}+\lambda_ig_i^*C' _i)$-numerically trivial, where $C'_{i}$ is the pullback of $C_{i}$. 
Therefore this sequence is a $(K_{Y_0}+\Delta_{Y_0})$-MMP over $Z'$ with scaling of $g_0^*C' _0$. 
Hence, this terminates by the case treated above. 
\end{proof}

\subsection{MMP}\label{ss3-lc-mmp}

The purpose of this subsection is to show Theorem \ref{t-lc-mmp}. 
We first treat the minimal model program for the pseudo-effective log canonical pairs 
(Theorem \ref{t-pseff-lc-mmp}). 
In this case, any $(K_X+\Delta)$-MMP with scaling of an ample divisor terminates. 

\begin{thm}\label{t-pseff-lc-mmp}
Let $k$ be a perfect field of characteristic $p>5$. 
Let $(X, \Delta)$ be a three-dimensional geometrically integral log pair over $k$. 
Let $f:X \to Z$ be a projective $k$-morphism to a quasi-projective $k$-scheme $Z$. 
Assume that there exist an $\R$-Cartier $\R$-divisor $C$ on $X$ and 
an effective $\R$-Cartier $\R$-divisor $C'$ on $X_{\overline k}:=X \times_k \overline k$ such that if $\Delta_{\overline k}$ and $C_{\overline k}$ denote the pullbacks of $\Delta$ and $C$ to $X_{\overline k}$ respectively, then 
\begin{enumerate}
\item[(a)] 
$K_X+\Delta$ is $f$-pseudo-effective, 
\item[(b)] 
$(X_{\overline k}, \Delta_{\overline k}+C')$ is log canonical, 
\item[(c)] 
$K_X+\Delta+C$ is $f$-nef, 
\item[(d)] 
$C$ is $f$-big, and 
\item[(e)] 
$C_{\overline k} \sim_{Z, \R} C'$. 
\end{enumerate}
Then the following hold. 
\begin{enumerate}
\item 
There exists no infinite sequence that is a $(K_X+\Delta)$-MMP over $Z$ with scaling of $C$. 
\item 
There exists a $(K_X+\Delta)$-MMP over $Z$ with scaling of $C$ that terminates. 
\end{enumerate}
\end{thm}

\begin{proof}
By Theorem \ref{t-mmp-scaling}, 
(1) implies (2). 
Hence it suffices to prove (1).

Assume that there exists an infinite sequence that is a $(K_X+\Delta)$-MMP over $Z$ with scaling of $C$: 
\begin{equation}\label{e-pseff-lc-mmp}
X=X_0 \dashrightarrow X_1 \dashrightarrow \cdots.
\end{equation}
Let us derive a contradiction. 
For the scaling coefficients $\lambda_0, \lambda_1, \ldots$, we set 
$\lambda:=\lim_{i\to \infty} \lambda_i$. 
If $\lambda \neq \lambda_i$ for any $i$, 
then the sequence (\ref{e-pseff-lc-mmp}) terminates 
by Corollary \ref{c-rel-lmm}. 
Hence, we may assume that $\lambda=\lambda_i$ for some $i$. 
Then the infinite sequence (\ref{e-pseff-lc-mmp}) is a $(K_X+\Delta+\frac{\lambda }{2} C)$-MMP over $Z$. 
Let 
\[
X_{\overline k}=X_{0, \overline k} \dashrightarrow X_{1, \overline k} \dashrightarrow \cdots 
\]
be the infinite sequence obtained by applying the base change $(-) \times_k \overline k$ 
to (\ref{e-pseff-lc-mmp}). 
Fix a projective birational morphism $g:Y^{(0)} \to X_{\overline k}=X_{0, \overline k}$ which satisfies the properties listed in 
Corollary \ref{c-dlt-modif}. 
In particular, if $\Delta_{Y^{(0)}}$ is the $\R$-divisor defined by $K_{Y^{(0)}}+\Delta_{Y^{(0)}}=g^*(K_{X_{\overline k}}+\Delta_{\overline k})$, then 
$(Y^{(0)}, \Delta_{Y^{(0)}})$ is a $\Q$-factorial dlt pair. 
Set $C_{Y^{(0)}}$ to be the pullback of $C'$ to $Y^{(0)}$. 
By standard argument, 
we can construct an infinite sequence that is a $(K_{Y^{(0)}}+\Delta_{Y^{(0)}}+\frac{1}{2} C_{Y^{(0)}})$-MMP over $Z$. 
As $K_X+\Delta+\frac{\lambda }{2} C$ is big over $Z$, 
$K_{Y^{(0)}}+\Delta_{Y^{(0)}}+\frac{1}{2} C_{Y^{(0)}}$ is big over $Z_{\overline k}$. 
This contradicts Theorem \ref{t-lc-eff-mmp2}. 
Therefore, (1) holds. 
\end{proof}

\begin{thm}\label{t-lc-mmp}
Let $k$ be a perfect field of characteristic $p>5$. 
Let $(X, \Delta)$ be a three-dimensional log canonical pair over $k$. 
Let $f:X \to Z$ be a projective $k$-morphism to a quasi-projective $k$-scheme $Z$. 
Then there exists a $(K_X+\Delta)$-MMP over $Z$ that terminates. 
\end{thm}

\begin{proof}
By standard argument, we may assume that $f_*\MO_X=\MO_Z$ and 
$k$ is algebraically closed in $K(X)$. 
In particular, $X$ is geometrically integral over $k$. 
Fix an algebraic closure $\overline k$ of $k$. 
Let $f_{\overline k}:X_{\overline k} \to Z_{\overline k}$ be the base change $f \times_k \overline k$.

If $K_X+\Delta$ is $f$-pseudo-effective, then the assertion follows from Theorem \ref{t-pseff-lc-mmp}. 
Thus we may assume that $K_X+\Delta$ is not $f$-pseudo-effective. 
Fix a projective birational morphism $g:Y \to X$ which satisfies the properties (1) and (2) of Corollary \ref{c-dlt-modif}. 
We prove the assertion of Theorem \ref{t-lc-mmp} by induction on $\rho(Y/Z)$. 

If $\rho(X/Z) \leq 1$, then the assertion holds. 
Indeed, if $\rho(X/Z)=0$, then there is nothing to show. 
If $\rho(X/Z)=1$, then the Stein factorisation of $f:X \to Z$ 
is a $(K_X+\Delta)$-Mori fibre space over $Z$. 

Therefore, we may assume that $\rho(X/Z)>1$. 
Pick an $f$-ample $\R$-Cartier $\R$-divisor $C$ on $X$ and 
an effective $\R$-Cartier $\R$-divisor $\overline C$ on $X_{\overline k}$ such that 
\begin{enumerate}
\item[(i)] 
$(X_{\overline k},\Delta_{\overline k}+\overline C)$ is log canonical, $K_{X}+\Delta+C$ is $f$-nef, and 
\item[(ii)] 
if we set
\[
\nu=\inf\{\mu \in \R_{\geq 0} \mid K_{X}+\Delta+\mu C \;{\rm is}\;f\mathchar`-{\rm pseudo}\mathchar`-{\rm effevtive}\},
\]
then $K_{X}+\Delta+\nu C \not \equiv_{f} 0$. 
\item[(iii)] 
$\overline C \sim_{\R, Z_{\overline k}} C_{\overline k}$, 
where $C_{\overline k}$ denotes the pullback of $C$ to $X_{\overline k}$. 
\end{enumerate}
Note that $0<\nu \leq 1$. 
By Theorem \ref{t-pseff-lc-mmp}, 
there exists a $(K_{X}+\Delta+\nu C)$-MMP over $Z$ with scaling of $C$ 
that terminates: 
\begin{equation}\label{e1-lc-mmp}
X=X_0 \dashrightarrow X_1 \dashrightarrow \cdots \dashrightarrow X_{\ell}=:X'. 
\end{equation}
If $\Delta'$ and $C'$ denote the proper transforms of $\Delta$ and $C$ on $X'$ respectively, 
then $K_{X'}+\Delta'+\nu C'$ is nef over $Z$. 
Moreover, by Theorem \ref{t-R-rel-bpf1}, we see that $K_{X'}+\Delta'+\nu C'$ is semi-ample over $Z$. 
We obtain projective morphisms: 
\[
X' \xrightarrow{f'} Z' \xrightarrow{\psi} Z, 
\]
where $f'$ is the projective morphism over $Z$ with $f'_*\MO_{X'}=\MO_{Z'}$ that is induced by $K_{X'}+\Delta'+\nu C'$, 
and $\psi$ is the induced morphism.  
By the choice of $\nu$, it holds that $K_{X'}+\Delta'+\nu C'$ is not big over $Z$. 
Therefore we get $\dim X' >\dim Z'$, which implies that 
$K_{X'}+\Delta'\equiv _{f'}-\nu C'$ is not $f'$-pseudo-effective. 
Recall that $g:Y \to X$ is a projective morphism which satisfies the properties (1) and (2) of Corollary \ref{c-dlt-modif}. 
Let $\Delta_Y$ be the $\R$-divisor defined by $K_{Y}+\Delta_Y=g^{*}(K_{X}+\Delta)$. 
Since the sequence (\ref{e1-lc-mmp}) is a $(K_{X}+\Delta)$-MMP over $Z$, 
we can construct $g':Y' \to X'$ and an $\mathbb{R}$-divisor $\Delta_{Y'}$ on $Y'$ such that 
\begin{itemize}
\item
$g'$ and $(Y',\Delta_{Y'})$ satisfies the properties (1) and (2) of Corollary \ref{c-dlt-modif}, and
\item
the induced birational map $\varphi:Y\dashrightarrow Y'$ is decomposed as a 
$(K_{Y}+\Delta_{Y})$-MMP over $Z$: 
\[
Y=Y^{(0)} \dashrightarrow \cdots \dashrightarrow Y^{(m)}=Y', 
\]
and $\varphi_{*}\Delta_Y=\Delta_{Y'}$.
\end{itemize}
In particular, we have $\rho(Y'/Z)\leq \rho(Y/Z)$. 

We prove that $\rho(Y'/Z')< \rho(Y/Z)$. 
When $K_{X'}+\Delta'+\nu C' \not \equiv_{\psi\circ f'}0$, 
there is a curve on $Z'$ whose image on $Z$ is a point. 
Therefore we get $\rho(Y'/Z')< \rho(Y'/Z)\leq \rho(Y/Z)$. 
When $K_{X'}+\Delta'+\nu C' \equiv_{\psi\circ f'}0$, the birational map $X\dashrightarrow X'$ contract a prime divisor $E$ because $K_{X}+\Delta+\nu C\not \equiv_{f}0$ by the above property (ii). 
Then 
\[
a_{E}(Y',\Delta_{Y'})=a_{E}(X',\Delta')>a_{E}(X,\Delta)=a_{E}(Y,\Delta_{Y}).
\]
Since $Y=Y^{(0)} \dashrightarrow \cdots \dashrightarrow Y^{(m)}=Y'$ is a $(K_{Y}+\Delta_{Y})$-MMP over $Z$, it contracts $g_{*}^{-1}E$. 
Therefore we have $\rho(Y'/Z')\leq \rho(Y'/Z)< \rho(Y/Z)$. 
In any case, we obtain $\rho(Y'/Z')< \rho(Y/Z)$. 

We apply the induction hypothesis to $(X',\Delta')$, $f':X' \to Z'$ and $g':Y' \to X'$. 
Since $K_{X'}+\Delta'$ is not $f'$-pseudo-effective, there is a $(K_{X'}+\Delta')$-MMP  over $Z'$ that terminates: 
\[
(X',\Delta')=(X_{\ell}, \Delta_{\ell}) \dashrightarrow (X_{\ell+1}, \Delta_{\ell+1}) \dashrightarrow \cdots \dashrightarrow (X_n, \Delta_n)=(X'',\Delta'').
\]  
Then there is a $(K_{X''}+\Delta'')$-Mori fibre space $g'':X''\to Z''$ over $Z'$. 
In particular, $g''$ is a $(K_{X''}+\Delta'')$-Mori fibre space over $Z$. 
Hence, the sequence 
\[
(X,\Delta)=(X_0, \Delta_0) \dashrightarrow \cdots \dashrightarrow (X_n, \Delta_n)=(X'',\Delta'')
\]
is a $(K_{X}+\Delta)$-MMP over $Z$ that terminates. 
\end{proof}

%%%%%%%%%%%%%%%%%%%%%%%%%%%%%%%%%%%%%%%%%%%%%%%%%%%%%%%%%%%%%%%%%%%%

\end{document}